\documentclass[final,leqno]{siamltex}

\usepackage{amssymb}
\usepackage{epsfig}
\usepackage{graphicx}
\usepackage{amsmath}
\usepackage{colortbl}
\usepackage{latexsym}
\usepackage{CJK}

\newtheorem{assumption}{Assumption}[section]
\newtheorem{problem}{Problem}[section]
\newtheorem{remark}{Remark}[section]

\newcommand{\qed}{\hbox{\rule[-2pt]{3pt}{6pt}}}

\allowdisplaybreaks

\title{Optimal Control and Stabilization for Linear Continuous-time Mean-field Systems with Delay}

\author{Xiao Ma$^\dagger$,~Qingyuan Qi$^\ddagger$, ~Xun Li$^\S$, and~Huanshui Zhang$^\dagger$\thanks{$^\dagger$School of Control Science and Engineering,
Shandong University, Jinan, Shandong, P.R. China 250061. H. Zhang is the
corresponding author (hszhang@sdu.edu.cn).
$^\ddag$Institute of Complexity Science, College of Automation, Qingdao University, Qingdao 266071, P.R. China (qiqy123@163.com).
$^\S$Department of Applied Mathematics, The Hong Kong Polytechnic University, Kowloon, Hong Kong, P.R. China (malixun@polyu.edu.hk).
}}

\begin{document}
\maketitle

\begin{abstract}
This paper studies optimal control and stabilization problems for continuous-time mean-field systems with input delay,
which are the fundamental development of control and stabilization problems for mean-field systems. There are two main contributions: 1) To the best of our knowledge, the present paper is first to establish the necessary and sufficient solvability condition for this kind of optimal control problem with delay, and to derive an optimal controller through overcoming the obstacle that separation principle no longer holds for multiplicative-noise systems; 2)  For the stabilization problem, under the assumption of exact observability,  we strictly prove that the system is stabilizable if and only if the algebraic Riccati equation has a unique positive definite solution.
\end{abstract}



\pagestyle{myheadings}
\thispagestyle{plain}
\markboth{TEX PRODUCTION}{USING SIAM'S \LaTeX\ MACROS}

\section{Introduction}

In this paper, we mainly investigate the optimal control and stabilization problem for linear continuous-time stochastic mean-field systems with time-delay. Time-delay systems have been extensively studied since 1950s because of its great importance in applications such as engineering and finance.
Lots of real-world systems evolve in accordance with not only their current state but also their previous history. It is a natural phenomenon that the past path-dependence exists in the fields such as chemistry, biology, finance, physics and economics. Therefore, the study for mean-field systems with delay is important yet meaningful in both theory and practice.
Pontryagin \cite{Pontryagin} initially established the maximum principle for control problems, Chyung-Lee \cite{CL} considered systems by linear differential equations involving time delays and obtained a variety of integral type criterion for optimality,
Artstein \cite{Artstein} studied systems with multiple delays, demonstrating a strong tool -- the reduction technique for stabilization problems by transforming linear systems with input delay into systems without delays.
Fridman-Shaked \cite{FS} and Park \cite{Park}
explored the stability criterion for systems with time-varying delays and time invariant uncertain delays.
Gozzi-Marinelli \cite{gozzi} introduced optimal advertising problem under uncertainty for the introduction of a new product to the market and adopted dynamic programming principle to construct optimal feedback controls.
In addition, the linear quadratic problems under stochastic systems with delay have been widely studied by a large number of research works. In particular, Zhang-Xu \cite{ZX} presented a complete solution for its optimal regulation and stabilization.
However, it is very difficult to tackle the optimal control problem for delayed systems owing to the infinite-dimensional state space structure.

The mean-field theory has been widely used in many fields, such as financial mathematics and statistical mechanics.
Stochastic mean-field systems involve state processes as well as their expected values.
The study of mean-field stochastic differential equations (MF-SDEs) goes back to the Mckean-Vlasov SDE proposed by Kac \cite{Kac} and McKean \cite{McKean} in the 1950s, which is the necessary condition of the optimal controller for continuous-time mean-field systems.
Since then, great efforts have been making on the exploration of MF-SDEs:
Buckdahn-Djehiche-Li-Peng \cite{BDLP} investigated a special mean-field problem in a purely stochastic approach,
Buckdahn-Li-Peng \cite{BLP}  obtained mean-field backward stochastic differential equations and related partial differential equations.
In the recent years, Yong \cite{Yong} paved the road for studying mean-field linear quadratic (LQ) controls, and presented optimal conditions and solutions to mean-field forward-backward stochastic equations (MF-FBSDEs) and Riccati equations for stochastic mean-field systems with deterministic coefficients, and Huang-Li-Yong \cite{HLY} established the stabilization criterion for mean-field LQ over an infinite horizon, derived optimal feedback controls via solutions to algebraic Riccati equations and further introduced a semi-definite programming method to fully tackle algebraic Riccati equations.

The problem that we consider in this paper is a continuation of Qi-Zhang \cite{QZ}, who studied problems of the delay-free continuous-time stochastic mean-field systems.  Qi-Zhang \cite{QZ} derived the necessary and sufficient solvability condition of the optimal control problem by solving the FBSDE and combining with the maximum principle, and achieved the stabilization condition for mean-field systems by an Lyapunov functional and the solution to the FBSDE.
Our study of the mean-field systems with input delay significantly differs from Qi-Zhang \cite{QZ} in the following aspects: First, the linear quadratic regulation is an important tool and one of the core concerns of modern control.
In the previous works, however, the linear quadratic regulation and stabilization problems for continuous-time Mean-field systems with input delay has not been studied yet. Second, a reduction method is applied to reduce the delayed system into a controlled system without delays, whose stabilizability is equivalent to the original system.

It should be emphasized that our work is a challenging, yet practically relevant and important, decision-making model with delay.
There exist great mathematical challenges, since it is
a fundamental development from delay-free mean-field systems to mean-field systems with delay. To be specific,  for the finite
horizon problem, it is very hard to obtain the solution to D-FBSDEs,
which is the essential for solving the control problem,  owing that the solution to the delayed
D-MFSDEs is nonlinear and the corresponding optimal controller involves
the conditional expectation of the state instead of the state itself.
It can be tackled by discretization method: Converting the system into a discrete
one, deriving the solution to the converted D-FBSDE, then taking limitation
of discrete-time solution yields the continuous-time solution.  In particular, for the infinite
horizon problem, the convergent analysis of the corresponding solution to the D-FBSEDs is quite complicated.

In addition, in the previous works, the control weighting matrix $R$ is assumed to be
positive definite when one seeks the stabilization conditions of infinite horizon including the standard LQ case. However, in this work, we  just assume that $R$ is positive semi-definite. The skills adopted in our paper are innovative, and the results established are significantly different from the existing theoretical aspect.

\newcommand{\RNum}[1]{\uppercase\expandafter{\romannumeral #1\relax}}

This paper is organized as follows. In Section \RNum{2}, for the finite horizon problem, the solvability condition of the optimal control and the analytical optimal controller are derived. In Section \RNum{3}, for the infinite horizon problem, the necessary and sufficient stabilization condition and the optimal stabilizing controller are presented. Section \RNum{4} provides a numerical example to shed light on the theoretical results established. We conclude the paper in Section \RNum{5}. Finally, the related proofs are presented in Appendix.

The following notation will be used throughout this paper.

\textbf{Notation}:
The equality between the stochastic variables and the uniqueness of the solution to the stochastic systems is in the sense of almost sure (a.s.);
$\mathbb{R}^n$ denotes the set of $n$-dimensional vectors;
$x'$ represents the transpose of $x$;
real matrix $A>0~(\geq 0)$ indicates that $A$ is positive definite (positive semi-definite);
$A^{-1}$ denotes the inverse of matrix $A$, and $A^\dag$ represents the Moore-Penrose inverse of $A$;
$\mathcal{F}_t$ denotes the filtration generated by the standard Brownian motion $\{w(s)|0\leq s\leq t\}$ and system initial state augmented by all the $\mathcal{P}$-null sets;
$\mathbb{E}[x(t)]$ means the unconditional expectation of $x(t)$;
$\{\Omega,~\mathcal{F},~\mathcal{P},~\{\mathcal{F}_t\}_{t\geq 0}\}$ is a complete stochastic probability space with natural filtration $\{\mathcal{F}_t\}_{t\geq 0}$;
$\mathbb{E}[\cdot|\mathcal{F}_t]$ denotes the conditional expectation with respect to $\mathcal{F}_t$;
$\frac{\partial f(x,y,\cdots)}{\partial x}$ denotes the partial derivative of $f$ with respect to $x$;
$\mathcal{I}_V$ is the indicator function of set $V$ satisfying
$\omega \in V,~\mathcal{I}_V = 1$, otherwise $\mathcal{I}_V = 0$;
$\langle\cdot,\cdot\rangle$ means the inner product in Hilbert space;
$\mathcal{M}_2$ denotes the Hilbert space $R^n\times L_2[-h,0]$ where $L_2[-h,0]$ is a square-integral function space;
$W^{1,2}(-h,0;R^m)$ means the Sobolev space of $R^m$-valued, absolutely continuous functions with square integrable derivatives on $[-h,0]$.

%

\section{Finite Horizon Problem}
\subsection{Problem Statement}
Let us consider the following continuous-time mean-field system with delay
\begin{align}\label{01}
\left\{\begin{array}{rl}
dx(t)=\!\!\! & \big\{Ax(t)+\bar{A}\mathbb{E}[x(t)]+Bu(t-h)+\bar{B}\mathbb{E}[u(t-h)]\big\}dt \\
&+\big\{Cx(t)+\bar{C}\mathbb{E}[x(t)]+Du(t-h)+\bar{D}\mathbb{E}[u(t-h)]\big\}dw(t), \\
x(0)=\!\!\! & x_0,~u(\tau)=\mu(\tau), \quad\quad \tau \in [-h,0),
\end{array}\right. 
\end{align}
where $x(t)$ is state process, $u(t)\in \mathbb{R}^m$ is control input process, $h ~ (>0)$ denotes input delay, $w(t)$ is one-dimensional standard Brownian motion. The initial values of the state and the input are given as $x_0$ and $\mu(\tau),~\tau\in[-h,0)$ respectively, where $\mu(\tau)$ is continuous and bounded over $[-h,0)$. In addition, $A,~\bar{A},~C,~\bar{C},~B,~\bar{B},~D,~\bar{D}$ are constant matrices with appropriate dimensions.

The corresponding mean-field type cost functional is given as
\begin{align}\label{02}
J_T= \mathbb{E}\bigg\{&\int_0^T x'(t)Qx(t)+\mathbb{E}[x'(t)]\bar{Q}\mathbb{E}[x(t)]dt \\ \notag
&+\int_h^T u'(t-h)Ru(t-h) +\mathbb{E}[u'(t-h)]\bar{R}\mathbb{E}[u(t-h)]dt \\ \notag
& +x'(T)P(T)x(T)+\mathbb{E}[x'(T)]\bar{P}(T)\mathbb{E}[x(T)]\bigg\},
\end{align}
where $Q,~\bar{Q},~R,~\bar{R},~P(T),~\bar{P}(T)$ are deterministic matrices with compatible dimensions.

\begin{assumption}\label{a01}
$Q,~Q+\bar{Q}\geq 0,~R,~R+\bar{R}\geq 0,~P(T)\geq 0,~P(T)+\bar{P}(T)\geq 0$.
\end{assumption}

The following section is discussed under Assumption \ref{a01}. The main problem over a finite time horizon can be stated as below:

\begin{problem}\label{p01}
Find an $\mathcal{F}_t$-adapted optimal controller of system \eqref{01} to minimize the cost functional \eqref{02}.
\end{problem}

\begin{remark}\label{r01}
The fundamental difficulty in solving Problem \ref{p01} is that separation principle no longer holds for multiplicative-noise systems.
\end{remark}

\subsection{Maximum Principle}
In order to solve Problem \ref{p01}, we introduce the maximum principle, which is the necessary condition for the optimal controller to minimize the cost functional.
\begin{theorem}\label{t01} \sl
If $u(t)$ is an $\mathcal{F}_t$-adapted optimal controller to Problem \ref{p01}, then it satisfies the following equilibrium equation:
\begin{align}\label{03}
&Ru(t)+\bar{R}\mathbb{E}[u(t)]+\mathbb{E}\Big[B'p(t+h)+\bar{B}'\mathbb{E}[p(t+h)] \\ \notag
&+D'q(t+h)+\bar{D}'\mathbb{E}[q(t+h)]\Big|\mathcal{F}_t\Big]=0,
\end{align}
where $p(t)$ and $q(t)$ satisfy the backward stochastic differential equation (BSDE) as follows:
\begin{align}\label{04}
\left\{\begin{array}{lll}dp(t)=-\{A'p(t)+\bar{A}'\mathbb{E}[p(t)]+C'q(t)+\bar{C}'\mathbb{E}[q(t)]+Qx(t)+\bar{Q}\mathbb{E}[x(t)]\}dt\\
~~~~~~~~~+q(t)dw(t), \\ [2mm]
p(T)=P(T)x(T)+\bar{P}(T)\mathbb{E}[x(T)].
\end{array}\right.
\end{align}
\end{theorem}

\textit{Proof}.
See Appendix A.
\qed

\subsection{Controller Design}
We introduce the definition of $\Pi(t,s)$ and some properties of derivations in the beginning of this section, in order to make the proof more explicit.
Set
\begin{align}\label{12}
&\Pi(t,s)=e^{A'(s-t)}\Pi(s,s)e^{A(s-t)}, \quad s\in[t,t+h], \\
&\Pi(t,t)=M'_1(t)[\Upsilon_1(t)]^{-1}M_1(t), \\
&\bar{\Pi}(t,t)=M'_2(t)[\Upsilon_2(t)]^{-1}M_2(t),
\end{align}
with the terminal values $\Pi(T,T+\theta)=\bar{\Pi}(T,T+\theta)=0~(\theta \geq 0) $, and the concrete forms of $M_i(t),~\Upsilon_i(t),i=1,2$ are determined in Theorem \ref{t02}.

In order to solve the MF-FBSED, $\Theta(t)$ and $d\Theta(t)$ are defined in Remark \ref{r1}.  Also, Lemma \ref{l1} is introduced to deal with the calculation of $d\Theta(t)$. 

\begin{lemma}\label{l1} \sl
Assumption \ref{a01} holds. Then, for $\hat{x}(t|s)=\mathbb{E}\{x(t)|\mathcal{F}_s\}$, we have 
\begin{align}\label{13}
\partial\hat{x}(t|s)=\big\{A\hat{x}(t|s)+\bar{A}\mathbb{E}[x(t)]+Bu(t-h)+\bar{B}\mathbb{E}[u(t-h)]\big\}dt.
\end{align}
\end{lemma}
\textit{Proof}.
See Appendix B.
\qed

\begin{remark}\label{r1}
For convenience, denote $\Theta(t)\doteq -\int_0^h \Pi(t,t+\theta)\hat{x}(t|t-h+\theta)d\theta$, it follows from Lemma \ref{l1} that we have
\begin{align}\label{18}
d\Theta(t)=&-d\bigg\{\int_t^{t+h}\Pi(t,\theta)\hat{x}(t|\theta-h)d\theta\bigg\} \\ \notag
=&\bigg\{\Pi(t,t)\hat{x}(t|t-h)-\Pi(t,t+h)x(t) \\ \notag
&+\int_t^{t+h}(A'\Pi(t,\theta)+\Pi(t,\theta)A)\hat{x}(t|\theta-h)d\theta \\ \notag
&-\int_t^{t+h}\Pi(t,\theta)\big[A\hat{x}(t|\theta-h)+\bar{A}\mathbb{E}[x(t)]+Bu(t-h)+\bar{B}\mathbb{E}[u(t-h)]\big]d\theta\bigg\}dt \\ \notag
=&\bigg\{\Pi(t,t)\hat{x}(t|t-h)-\Pi(t,t+h)x(t)-A'\Theta(t)\\\notag
&-\int_t^{t+h}\Pi(t,\theta)\big[\bar{A}\mathbb{E}[x(t)]+Bu(t-h)+\bar{B}\mathbb{E}[u(t-h)]\big]d\theta\bigg\}dt\\\notag
=&\Theta_1(t)dt. \\ \notag
\end{align}
\end{remark}

The solvability condition of Problem \ref{p01} can be presented in the following.

\begin{theorem}\label{t02} \sl
Assumption \ref{a01} holds. Then Problem \ref{p01} is uniquely solvable if and only if $\Upsilon_1(t)>0$ and $\Upsilon_2(t)>0$ for $t\in[h,T]$, where $\Upsilon_1(t)$ and $\Upsilon_2(t)$ satisfy
\begin{align}\label{05}
&\Upsilon_1(t)=R+D'P(t)D, \\
&\Upsilon_2(t)=R+\bar{R}+(D+\bar{D})'P(t)(D+\bar{D}),
\end{align}
and $(P(t), \bar{P}(t))$ is the solution of the following coupled Riccati equations:
\begin{align}\label{06}
-\dot{P}(t)= ~ & Q+P(t)A+A'P(t)+C'P(t)C-\Pi(t,t+h), \\
-\dot{\bar{P}}(t)= ~ & \bar{Q}+P(t)\bar{A}+\bar{A}'P(t)+\bar{P}(t)(A+\bar{A})+(A+\bar{A})'\bar{P}(t)+C'P(t)\bar{C}\\ \notag
& +\bar{C}'P(t)C +\bar{C}'P(t)\bar{C}-\bar{A}'\int_0^h\Pi(t,t+\theta)d\theta-\int_0^h\Pi(t,t+\theta)d\theta\bar{A}\\\notag
&+\Pi(t,t)-\bar{\Pi}(t,t), \end{align}
with the terminal conditions as those in \eqref{02}, and $M_1(t)$ and $M_2(t)$ are presented by
\begin{align}\label{07}
M_1(t)=&B'P(t)-B'\int_0^h\Pi(t,t+\theta)d\theta+D'P(t)C \quad \mbox{and} \\
M_2(t)=&(B+\bar{B})'(P(t)+\bar{P}(t))-(B+\bar{B})'\int_0^h\Pi(t,t+\theta)d\theta\\\notag
&+(D+\bar{D})'P(t)(C+\bar{C}),
\end{align}
respectively.
Then the optimal controller can be presented by
\begin{equation}\label{08}
u(t-h)=K(t)\hat{x}(t|t-h)+\bar{K}(t)\mathbb{E}[x(t)],
\end{equation}
where
\begin{align}\label{09}
\left\{\begin{array}{l}
K(t)=-[\Upsilon_1(t)]^{-1}M_1(t), \\
\bar{K}(t)=-\Big\{[\Upsilon_2(t)]^{-1}M_2(t)-[\Upsilon_1(t)]^{-1}M_1(t)\Big\},
\end{array}\right.
\end{align}
and the optimal cost functional can be represented by
\begin{align}\label{091}
J_T^*=&\mathbb{E}\Big\{\int_0^h\bigg[x'(t)Qx(t)+\mathbb{E}[x'(t)]\bar{Q}\mathbb{E}[x(t)]\Big]dt \\ \notag
& \quad +x'(h)\Big[P(h)x(h)+\bar{P}(h)\mathbb{E}[x(h)]-\int_0^h\Pi(h,h+\theta)\hat{x}(h|\theta)d\theta\Big]\bigg\}.
\end{align}
Meanwhile, the costate $p(t)$ and the state $x(t)$ satisfy the following relation:
\begin{align}\label{10}
p(t)=P(t)x(t)+\bar{P}(t)\mathbb{E}[x(t)]-\int_0^h\Pi(t,t+\theta)\hat{x}(t|t+\theta-h)d\theta.
\end{align}
\end{theorem}

\textit{Proof}.
See Appendix C.
\qed

\begin{remark} \label{r03}
It should be emphasized that in Theorem \ref{t01} we consider the case with $R\geq 0$ and $R+\bar{R}\geq 0$, which is a more general assumption than the requirement of positive definiteness. The results established in Theorem \ref{t01} can be degenerated to a special case if we let $h=0$ (traditional mean-field control) or $\bar{A}=\bar{B}=\bar{C}=\bar{D}=0$ (delayed systems control, see reference \cite{ZX}).
\end{remark}

\section{Infinite Horizon Problem}

\subsection{Problem Statement}

In this section, the optimal control and stabilization problems of infinite horizon will be studied. The cost functional of this problem is described as
\begin{align}\label{37}
J=\mathbb{E}\bigg\{\int_0^\infty x'(t)Qx(t)+\mathbb{E}[x'(t)]\bar{Q}\mathbb{E}[x(t)]+u'(t)Ru(t)+\mathbb{E}[u'(t)]\bar{R}\mathbb{E}[u(t)]dt\bigg\}.
\end{align}

We introduce some definitions and standard assumptions before stating the main problem.
\begin{definition}\label{d1}
The mean-field stochastic system is presented by
\begin{align}\label{59}
\left\{\begin{array}{l}
dx(t)=\big\{Ax(t)+\bar{A}\mathbb{E}[x(t)]\big\}dt+\big\{Cx(t)+\bar{C}\mathbb{E}[x(t)]\big\}dw(t), \\
Y_t=\mathcal{Q}^{\frac{1}{2}}X(t),
\end{array}\right. 
\end{align}
where
$$X(t)=\left(
  \begin{array}{ccc}
    x(t)-\mathbb{E}[x(t)] \\
    \mathbb{E}[x(t)] \\
 \end{array}
\right)
\quad \mbox{and} \quad
\mathcal{Q}=\left(
  \begin{array}{ccc}
    Q&0 \\
    0&Q+\bar{Q} \\
   \end{array}
\right).$$
System \eqref{59}, or $(A,~\bar{A},~C,~\bar{C},\mathcal{Q}^{\frac{1}{2}})$ for short, is said to be exactly observable, if for any $T>0$,
\begin{align}\label{60}
Y_t=0,~a.s.,~\forall t\in[0,T]~~~~\Rightarrow~~~~x_0=0.
\end{align}
\end{definition}

\begin{definition}\label{d2}
System \eqref{01} with $u=0$ is called asymptotically mean-square stable if $\displaystyle\lim_{t\rightarrow \infty}\mathbb{E}[x'(t)x(t)]=0$ holds for any initial values $x_0$.
\end{definition}

\begin{definition}\label{d3}
System \eqref{01} is called stabilizable in the mean-square sense if there exists an $\mathcal{F}_t$-adapted controller $u(t)$ such that the closed-loop system is asymptotically mean-square stable.
\end{definition}

\begin{assumption}\label{a2}
$Q,~Q+\bar{Q}\geq0,~R,R+\bar{R}\geq 0$.
\end{assumption}

\begin{assumption}\label{a3}
$(A,\bar{A},C,\bar{C},\mathcal{Q}^{\frac{1}{2}})$ is exactly observable.
\end{assumption}

The main problem over an infinite time horizon can be described as follows:
\begin{problem}\label{p02}
Find an $\mathcal{F}_t$-adapted controller $u(t)$ to minimize cost functional \eqref{37} and stabilize system \eqref{01} in the mean-square sense.
\end{problem}

\subsection{Solution to Problem \ref{p02}}

For convenience of discussion, in this section, we re-denote $P(t),~\bar{P}(t),~\Pi(t,t+\theta)$ and $\bar{\Pi}(t,t)$ as $P_T(t),~\bar{P}_T(t)$, $\Pi_T(t,t+\theta)$ and $\bar{\Pi}_T(t,t)$ to make the terminal time $T$ explicit. Also, the terminal values are set as $P_T(T)=\bar{P}_T(T)=\Pi_T(T,T+\theta)=\bar{\Pi}_T(T,T)=0$. Then we have the following theorem.

\begin{theorem}\label{t03} \sl
Under Assumptions \ref{a2} and \ref{a3}, if system \eqref{01} is mean-square stabilizable, then the following two assertions hold:

(1)~$P_T(t),~\bar{P}_T(t),~\Pi_T(t,t+\theta),~\theta\in[0,h]$ and $\bar{\Pi}_T(t,t)$ are convergent when $T\rightarrow \infty$, and their limits,
denoted by $P,~\bar{P}$, $\Pi(\theta)$ and $\bar{\Pi}(0)$, satisfy the following coupled algebraic Riccati equations:
\begin{align}\label{38}
-Q=&PA+A'P+C'PC-\Pi(h), \\\label{38b}
-\bar{Q}=&P\bar{A}+\bar{A}'P+\bar{P}(A+\bar{A})+(A+\bar{A})'\bar{P}+C'P\bar{C}+\bar{C}'PC+\bar{C}'P\bar{C} \\ \notag
&-\bar{A}'\int_0^h\Pi(\theta)d\theta-\int_0^h\Pi(\theta)d\theta\bar{A}
+\Pi(0)-\bar{\Pi}(0),
\end{align}
where
\begin{align}\label{39}
 &\Pi(\theta)=e^{A'\theta}\Pi(0)e^{A\theta}, \\ \label{39a}
 &\Pi(0)=M'_1\Upsilon_1^{\dag}M_1, \\ \label{39b}
&\bar{\Pi}(0)=M'_2\Upsilon_2^{\dag}M_2,
\end{align}
with the regular condition
\begin{align}\label{42}
\Upsilon_i\Upsilon_i^{\dag}M_i=M_i,~~~i=1,2,
\end{align}
and
\begin{align}\label{40}
&M_1=B'P-B'\int_0^h\Pi(\theta)d\theta+D'PC, \\ \label{40a}
&M_2=(B+\bar{B})'(P+\bar{P})-(B+\bar{B})'\int_0^h\Pi(\theta)d\theta+(D+\bar{D})'P(C+\bar{C}), \\ \label{40b}
&\Upsilon_1=R+D'PD, \\ \label{40c}
&\Upsilon_2=R+\bar{R}+(D+\bar{D})'P(D+\bar{D}).
\end{align}

(2)~ The following results hold
\begin{align}\label{41}
P-\int_0^h\Pi(\theta)d\theta>0,~~~~P+\bar{P}-\int_0^h\Pi(\theta)d\theta>0.
\end{align}
\end{theorem}

\textit{Proof}.
See Appendix D.
\qed

\begin{remark}
If the coupled algebraic Riccati equations \eqref{38}-\eqref{38b} have a solution $(P, \bar{P})$ satisfying
$$P-\int_0^h\Pi(\theta)d\theta>0, \quad P+\bar{P}-\int_0^h\Pi(\theta)d\theta>0,$$
then the solution is positive definite. It follows from Theorem \ref{t03} that if system \eqref{01} is mean-square stabilizable, then \eqref{38}-\eqref{38b} have a positive definite solution.
\end{remark}

\begin{lemma}\label{l2} \sl
Under Assumptions \ref{a2} and \ref{a3}, the following system $(\tilde{A},~\tilde{C},~\tilde{\mathcal{Q}}^{\frac{1}{2}})$ is exactly observable:
\begin{align}\label{76}
\left\{\begin{array}{l}
dX(t)=\tilde{A}X(t)dt+\tilde{C}X(t)dw(t), \\
\tilde{\mathcal{Y}}(t)=\tilde{\mathcal{Q}}^{\frac{1}{2}}X(t),
\end{array}\right.
\end{align}
where $X(t)=\left(
  \begin{array}{ccc}
    x(t)-\mathbb{E}[x(t)] \\
    \mathbb{E}[x(t)] \\
 \end{array}
\right)$, for $t\in[0, T]$,

$\tilde{A}=\left[
  \begin{array}{ccc}
    A+BK&0 \\
    0&A+\bar{A}+(B+\bar{B})(K+\bar{K}) \\
  \end{array}
\right]
$,

$\tilde{C}=\left[
  \begin{array}{ccc}
    C+DK&C+\bar{C}+(D+\bar{D})(K+\bar{K}) \\
    0&0 \\
  \end{array}
\right]
$,

$\tilde{\mathcal{Q}}=\left[
  \begin{array}{ccc}
    Q+K'RK&0 \\
    0&Q+\bar{Q}+(K+\bar{K})'(Q+\bar{Q})(K+\bar{K}) \\
  \end{array}
\right]\geq 0
$,
i.e., for any $T \geq 0$, $\tilde{\mathcal{Y}}(t)=0,~a.s.,~\forall t\in [0,~T]\Rightarrow x_0=0$.
\end{lemma}

\textit{Proof}.
See \cite{QZ}.
\qed

We have the following main result.

\begin{theorem}\label{t04} \sl
Under Assumptions \ref{a2} and \ref{a3}, system \eqref{01} is stabilizable in the mean-square sense if and only if the coupled algebraic Riccati equations \eqref{38}-\eqref{38b} have a unique positive definite solution $(P, \bar{P})$ such that $P-\int_0^h\Pi(\theta)d\theta > 0$ and $P+\bar{P}-\int_0^h\Pi(\theta)d\theta > 0$.
In this case, the optimal stabilizing controller is given by
\begin{align}\label{59b}
u(t-h)=K\hat{x}(t|t-h)+\bar{K}\mathbb{E}[x(t)],
\end{align}
where $K$ and $\bar{K}$ are given by
\begin{align}\label{60b}
&K=-\Upsilon_1^{\dag}M_1,\\
&\bar{K}=-\{\Upsilon_2^{\dag}M_2-\Upsilon_1^{\dag}M_1\},
\end{align}
and $M_i,~\Upsilon_i,~i=1,2$ are given by \eqref{40}-\eqref{40c}. Then the corresponding optimal cost functional is presented by
\begin{align}\label{61}
J^*=&\mathbb{E}\bigg\{x'(0)\bigg[Px(0)+\bar{P}\mathbb{E}[x(0)]-\int_0^h\Pi(\theta)d\theta x(0)\bigg] \\ \notag
&+\int_0^h\bigg[-\mathbb{E}[x'(t)]\Pi(0)\mathbb{E}[x(t)]+\mathbb{E}[x'(t)]\bar{\Pi}(0)\mathbb{E}[x(t)]\\ \notag
&+\hat{x}'(t|t-h)\Pi(0)\hat{x}(t|t-h)+2x'(t)M'_1u(t-h) \\ \notag
&+2\mathbb{E}[x'(t)][M_2-M_1]'\mathbb{E}[u(t-h)]+u'(t-h)D'PDu(t-h) \\ \notag
&+\mathbb{E}[u'(t-h)][\bar{D}'PD+D'P\bar{D}+\bar{D}'P\bar{D}]\mathbb{E}[u(t-h)]\bigg]dt \bigg\}.
\end{align}
\end{theorem}

\textit{Proof}.
See Appendix E.
\qed

\begin{remark}\label{r00001}
It is worth noting that the stabilization results in Theorem \ref{t04} are derived under the condition of $R\geq 0,~R+\bar{R}\geq0$, which is a more general condition than the standard assumption of $R>0$ in classical control field.
\end{remark}

\section{Numerical Example}

\begin{figure}[htbp]
  \centering
  \includegraphics[width=0.48\textwidth]{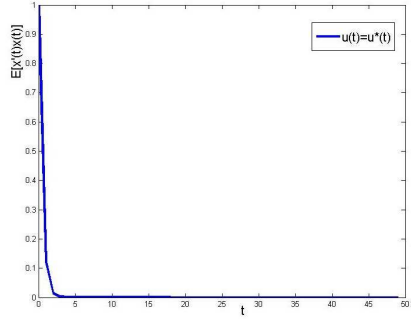}
  \caption{Simulations for the system state trajectory E[x'(t)x(t)].}\label{fig1}
\end{figure}

In this section, we verify the obtained results using the following numerical example.
Consider system \eqref{01} and cost functional \eqref{37} with
\begin{align*}
A=\bar{A}=2,~B=\bar{B}=1,~C=\bar{C}=0.2,~D=\bar{D}=0.1,~h=1,~Q=\bar{Q}=R=\bar{R}=1,
\end{align*}
and the initial values $x(0)=1,~u(\tau)=0,~\tau\in[-1,0)$.
It is obvious that these parameters satisfy Assumption \ref{a3}.
From \eqref{38}-\eqref{38b}, solving the coupled algebraic Riccati equation yields
$P\approx 1220.2,~\bar{P}\approx2.90$ and$\int_0^1\Pi(\theta)d\theta\approx1210$, therefore the unique solution satisfies $P-\int_0^1\Pi(\theta)d\theta\approx 10.2>0$ and $P+\bar{P}-\int_0^1\Pi(\theta)d\theta\approx13.1>0$. It follows from Theorem \ref{t04} that the system is stabilizable in the mean-square sense and the corresponding optimal stabilizing controller is calculated approximatively as $u(t-1)\approx-2.6211\hat{x}(t|t-1)+0.1842\mathbb{E}[x(t)]$. The simulation of numerical results is shown in Fig \ref{fig1}.

\section{Conclusion}

The LQ regulation and stabilization problems for continuous-time mean-field systems with delay is investigated in this paper. An analytical controller is derived by solving the D-MFSDE and using maximum principle. The necessary and sufficient stabilization condition can be proposed by an Lyapunov functional with the optimal cost functional. We prove that under the assumption of exact observability, the system is stabilizable if and only if the algorithm of Riccati equation admits a unique positive definite solution.

\section{Appendix}
\subsection{Appendix A: Proof of Theorem \ref{t01}}
\textit{Proof}.
Convex variation method will be applied in the following proof.
For any $u(t),~\delta u(t)\in \mathcal{U}[0,T-h]$, let $u^{\varepsilon}(t)=u(t)+\varepsilon \delta u(t),~\varepsilon\in(0,1)$. The corresponding state and cost functional with $u^{\varepsilon}(t)$ are denoted as $x^{\varepsilon}(t)$ and $J^{\varepsilon}_T$, respectively.
Let $X(t)=\left[
  \begin{array}{ccc}
    x(t) \\
    \mathbb{E}[x(t)] \\
  \end{array}
\right]
$,
then the variation of $X(t)$ satisfies
\begin{align}\label{r01b}
 &d\delta X(t)
 =d\left[
  \begin{array}{ccc}
    x^{\varepsilon}(t)-x(t) \\
    \mathbb{E}[x^{\varepsilon}(t)]-\mathbb{E}[x(t)] \\
  \end{array}\right]\\\notag
= & \Big[\mathcal{A}\delta X(t)+\mathcal{B}\varepsilon\delta u(t-h)+\mathcal{\bar{B}}\varepsilon\delta\mathbb{E}[u(t-h)]\Big]dt \\ \notag
&+ \Big[\mathcal{C}\delta X(t)+\mathcal{D}\varepsilon\delta u(t-h)+\mathcal{\bar{D}}\varepsilon\delta\mathbb{E}[u(t-h)]\Big]dw(t),
\end{align}
where
\[
\begin{array}{lll}
\mathcal{A}=\left[
  \begin{array}{ccc}
    A&\bar{A} \\
    0&A+\bar{A} \\
  \end{array}
\right],
& \mathcal{B}=\left[
  \begin{array}{ccc}
    B \\
    0 \\
  \end{array}
\right],
& \mathcal{\bar{B}}=\left[
  \begin{array}{ccc}
    B \\
    B+\bar{B} \\
  \end{array}
\right], \\ [5mm]
\mathcal{C}=\left[
  \begin{array}{ccc}
  C&\bar{C} \\
    0&0 \\
  \end{array}
\right],
& \mathcal{D}=\left[
  \begin{array}{ccc}
    D \\
    0 \\
  \end{array}
\right],
& \mathcal{\bar{D}}=\left[
  \begin{array}{ccc}
    \bar{D} \\
    0 \\
  \end{array}
\right].
\end{array}
\]
From Yong \cite{YZ} (Theorem 6.14), $\delta X(t)$ satisfies
\begin{align}\label{r02}
\delta x(t)=&\left[
  \begin{array}{ccc}
    I_n&0
  \end{array}
\right]\Phi(t)\int_h^t\Phi^{-1}(s)\Big[\mathcal{B}\varepsilon\delta u(s-h)+\mathcal{\bar{B}}\mathbb{E}[\varepsilon\delta u(s-h)]\\\notag
&-\mathcal{C}(\mathcal{D}\varepsilon\delta u(s-h)+\mathcal{\bar{D}}\mathbb{E}[u(s-h)])\Big]ds,
\end{align}
where $\Phi(t)$ is the unique solution of the following SDE
\begin{equation}\label{r03b}
 \left\{\begin{array}{lll}d\Phi(t)=\mathcal{A}\Phi(t)dt+\mathcal{C}\Phi(t)dw(t),\\\notag
\Phi(0)=I_{2n}, &\textrm{$$}\end{array}\right.
\end{equation}
and there exists $\Psi(t)=\Phi^{-1}(t)$ satisfying
\begin{equation}\label{r04}
\left\{\begin{array}{lll}d\Psi(t)=-\Psi(t)(\mathcal{A}-\mathcal{C}^2)dt-\Psi(t)Cdw(t), &\textrm{$$}\\\notag
\Psi(0)=I_{2n}. &\textrm{$$}\end{array}\right.
\end{equation}
The variation of $J_T$ satisfies:
\begin{align}\label{r05}
\delta J_T=&J_T^\varepsilon-J_T \\ \notag
=&2\mathbb{E}\bigg\{\int_0^T\Big[x'(t)Q+\mathbb{E}[x'(t)]\bar{Q}\Big]\delta x(t)dt \\ \notag
&+\int_h^T\Big[u'(t-h)R+\mathbb{E}[u'(t-h)]\bar{R}\Big]\varepsilon\delta u(t-h)dt \\ \notag
&+\Big[x'(T)P(T)+\mathbb{E}[x'(T)]\Big]\delta x(T)\bigg\}+o(\varepsilon),
\end{align}
where $o(\varepsilon)$ denotes infintesimal of high order with $\varepsilon$.

Substituting \eqref{r02} into \eqref{r05} yields
\begin{align}\label{r06}
&\begin{array}{rl}
\frac{1}{2}\delta J_T
=\!\!\!\!\!&\mathbb{E}\bigg\{\displaystyle\int_h^T\bigg(\int_s^T \Big(x'(t)Q+\mathbb{E}[x'(t)]\bar{Q}\Big)\left[
  \begin{array}{ccc}
    I_n&0\\
  \end{array}
\right]\Phi(t)dt \\ \notag
&+\Big(x'(T)P(T)+\mathbb{E}[x'(T)]\bar P(T)\Big)\left[
  \begin{array}{ccc}
    I_n&0\\
  \end{array}
\right] \Phi(T)\Phi^{-1}(s)(\mathcal{B}-\mathcal{CD}) \\ \notag
&+\mathbb{E}\bigg[\displaystyle\int_s^T \Big(x'(t)Q+\mathbb{E}[x'(t)]\bar{Q}\Big)\left[
  \begin{array}{ccc}
    I_n&0\\
  \end{array}
\right]\Phi(t)dt \\ \notag
&+
\Big(x'(T)P(T)+\mathbb{E}[x'(T)]\bar P(T)\Big)\left[
  \begin{array}{ccc}
    I_n&0\\
  \end{array}
\right]\Phi(T)]\Phi^{-1}(s)(\bar{B}-C\bar{D})\bigg] \\ \notag
&+u'(s-h)R+\mathbb{E}[u'(s-h)]\bar{R}\bigg)\varepsilon\delta u(s-h)ds\bigg\} \\ \notag
&+\mathbb{E}\bigg\{\bigg[\displaystyle\int_0^T \Big(x'(t)Q+\mathbb{E}[x'(t)]\bar{Q}\Big)\left[
  \begin{array}{ccc}
    I_n&0\\
  \end{array}
\right]\Phi(t)dt\\\notag
&+
\Big(x'(T)P(T)+\mathbb{E}[x'(T)]\bar P(T)\Big)
\left[
  \begin{array}{ccc}
    I_n&0\\
  \end{array}
\right]\Phi(T)\bigg] \\ \notag
&\times\displaystyle\int_h^T\Phi^{-1}(s)(D\varepsilon\delta u(s-h)+\bar{D}\mathbb{E}[\varepsilon\delta u(s-h))]dw(s)\bigg\}+o(\varepsilon).
\end{array}\\
\end{align}

Let $\xi=\int_0^T (x'(t)Q+\mathbb{E}[x'(t)]\bar{Q})\left[
  \begin{array}{ccc}
    I_n&0\\
  \end{array}
\right]\Phi(t)dt+
(x'(T)P(T)+\mathbb{E}[x'(T)])\left[
  \begin{array}{ccc}
    I_n&0\\
  \end{array}
\right]\Phi(T)$,
similar to the proof in \cite{QZ}, $\mathbb{E}[\xi|\mathcal{F}_s]$ is a martingale with respect to $s$, there exists a unique $\mathcal{F}_t$-adapted process $\eta(t)$ such that
\begin{align}\label{r07}
\xi=\mathbb{E}[\xi]+\int_0^T\eta'(t)dw(t),
\end{align}
and $\delta J_T$ can be rewitten as:
\begin{align}\label{r08}
\frac{1}{2}J_T=&\mathbb{E}\int_h^T H'(s)\varepsilon\delta u(s-h)ds+o(\varepsilon) \\ \notag
=&\mathbb{E}\int_h^T \mathbb{E}[H'(s)|\mathcal{F}_{s-h}]\varepsilon\delta u(s-h)ds+o(\varepsilon),
\end{align}
where
\begin{align}\label{r09}
\begin{array}{rl}
H'(s)=\!\!\!&\displaystyle\int_s^T (x'(t)Q+\mathbb{E}[x'(t)]\bar{Q})\left[
  \begin{array}{ccc}
    I_n&0\\
  \end{array}
\right]\Phi(t)dt \\ \notag
& +
\Big(x'(T)P(T)+\mathbb{E}[x'(T)]\bar P(T)\Big)
\left[
  \begin{array}{ccc}
    I_n&0\\
  \end{array}
\right]\Phi(T)]\Phi^{-1}(s)(\mathcal{B}-\mathcal{CD}) \\ \notag
&+\mathbb{E}\bigg\{\displaystyle\int_s^T (x'(t)Q+\mathbb{E}[x'(t)]\bar{Q})\left[
  \begin{array}{ccc}
    I_n&0 \\
  \end{array}
\right]\Phi(t)dt \\ \notag
&+\Big(x'(T)P(T)+\mathbb{E}[x'(T)]\Big)\left[
  \begin{array}{ccc}
    I_n&0 \\
  \end{array}
\right]\Phi(T)]\Phi^{-1}(s)(\mathcal{\bar{B}}-\mathcal{C\bar{D}})\bigg\} \\ \notag
&+u'(s-h)R+\mathbb{E}[u'(s-h)]\bar{R}+\eta'(s)\Phi^{-1}(s)\mathcal{D}+\mathbb{E}[\eta'(s)\Phi^{-1}(s)\mathcal{\bar{D}}].
\end{array} \\
\end{align}

Thereby, the necessary condition of minimizing $J_T$ is $\mathbb{E}[H'(s)|\mathcal{F}_{s-h}]=0$ due to the arbitrariness of $\delta u(s-h)$.

Similarly, by defining $(\check{p}(s),\check{q}(s))$ appropriately, the necessary condition can be rewritten as:
\begin{align}\label{r10}
Ru(s-h)+\bar{R}\mathbb{E}u(s-h)+\mathbb{E}\Big[\mathcal{B}'\check{p}(s)+\mathcal{D}'\check{q}(s)+\mathbb{E}[\mathcal{\bar{B}}'\check{p}(s)+\mathcal{\bar{D}}'\check{q}(s)]\Big\vert\mathcal{F}_{s-h}\Big] = 0,
\end{align}
and $(\check{p}(s),\check{q}(s))$ satisfies the following BSDE:
\begin{align}\label{R11}
\left\{\begin{array}{lll}d\check{p}(t)=-\bigg\{\left[
  \begin{array}{ccc}
    A&\bar{A} \\
    0&A+\bar{A}\\
  \end{array}
\right]'\check{p}(t)
+\left[
  \begin{array}{ccc}
    C&\bar{C} \\
    0&0 \\
  \end{array}
\right]'\check{q}(t)
+\left[
  \begin{array}{ccc}
    I_n \\
    0 \\
  \end{array}
\right]\big(Qx(t)+\bar{Q}\mathbb{E}[x(t)]\big)
\bigg\}dt \\ [5mm]
~~~~~~~~~~~+\check{q}(t)dw(t), &\textrm{$$} \\ [5mm]
\check{p}(T)=\left[
  \begin{array}{ccc}
    P(T)&\bar{P}(T) \\
    0&0 \\
  \end{array}
\right]\left[
  \begin{array}{ccc}
    x(T) \\
    \mathbb{E}[x(T)] \\
  \end{array}
\right]. &\textrm{$$}\end{array}\right.
\end{align}
According to \cite{RB} and \cite{RB2}, there exists a unique solution $[\check{p}(\cdot),~\check{q}(\cdot)]_{2n\times 1}=[(p(\cdot),p_2(\cdot))'$, $(q(\cdot),q_2(\cdot))']$ to \eqref{04}.
Let $[\hat{p}(\cdot),\hat{q}(\cdot)]=[(p (\cdot),\mathbb{E}[p(\cdot)])',(q(\cdot),\mathbb{E}[q(\cdot)])']$, it can be verified that $(\hat{p}(\cdot),\hat{q}(\cdot))$ satisfies \eqref{R11} according to the uniqueness of the solution. Hence, we can conclude that \eqref{04} holds and the optimal condition is as \eqref{03}. The proof is completed.
\qed

\subsection{Appendix B: Proof of Lemma \ref{l1}}
\textit{Proof}.
Using the variation of constants formula for linear SDEs, we express $x(t)$ in \eqref{01} as
\begin{align}\label{14}
x(t)=&\Phi(t)[\Phi(s)]^{-1}x(s)+\Phi(t)\int_s^t[\Phi(\theta)]^{-1}\Big[\bar{A}\mathbb{E}[x(\theta)]+Bu(\theta-h)\\ \notag
&+\bar{B}\mathbb{E}[u(\theta-h)]
-C(\bar{C}\mathbb{E}[x(\theta)]+Du(\theta-h)+\bar{D}\mathbb{E}[u(\theta-h)])\Big]d\theta\\\notag
&+\Phi(t)\int_s^t\Big[\Phi(\theta)]^{-1}[\bar{C}\mathbb{E}[x(\theta)] +Du(\theta-h)+\bar{D}\mathbb{E}[u(\theta-h)]\Big]dw(\theta),
\end{align}
where
\begin{align}\label{15}
\Phi(t)=e^{At}I+\int_0^te^{A(t-\theta)}C\Phi(\theta)dw(\theta).
\end{align}
Substituting \eqref{15} into \eqref{14} and taking the expectation yields
\begin{align}\label{16}
&\mathbb{E}\left\{\Phi(t)\int_s^t\Big[\Phi(\theta)]^{-1}[\bar{C}\mathbb{E}[x(\theta)]+Du(\theta-h)+\bar{D}\mathbb{E}[u(\theta-h)]\Big]dw(\theta)\Big\vert\mathcal{F}_s\right\} \\ \notag
=&\mathbb{E}\left\{\int_0^te^{A(t-\theta)}C\Phi(\theta)dw(\theta)\int_s^t[\Phi(\theta)]^{-1}\Big[\bar{C}\mathbb{E}[x(\theta)]+Du(\theta-h)\right. \\ \notag
&\left.+\bar{D}\mathbb{E}[u(\theta-h)]\Big]dw(\theta)\Big\vert\mathcal{F}_s\right\} \\ \notag
=&\mathbb{E}\left\{\int_s^te^{A(t-\theta)}C\Big[\bar{C}\mathbb{E}[x(\theta)]+Du(\theta-h)+\bar{D}\mathbb{E}[u(\theta-h)]\Big]d\theta\right\}.
\end{align}
Using \eqref{14} and \eqref{16}, we have
\begin{align}\label{17}
\hat{x}(t|s)=e^{A(t-s)}x(s)+\int_s^te^{A(t-\theta)}\Big[\bar{A}\mathbb{E}[x(\theta)]+Bu(\theta-h)+\bar{B}\mathbb{E}[u(\theta-h)]\Big]d\theta.
\end{align}
Then, taking partial derivation of $\hat{x}(t|s)$ with respect to $t$, we obtain the desired result \eqref{13}.

\subsection{Appendix C: Proof of Theorem \ref{t02}}

\textit{Proof}.
``Sufficiency":
Firstly, from the result of discrete-time and early work, we claim that the costate and state satisfy \eqref{10}. Next, we will prove that Problem \ref{p01} is uniquely solvable under the assumption of  $\Upsilon_1(t)>0$ and $\Upsilon_2(t)>0$, $t\in[h,T]$ hold.
In fact, applying It\^{o}'s formula to $p'(t)x(t)$, we get
\begin{align}\label{26}
& d[p'(t)x(t)] \\ \notag
=&[dp(t)]'x(t)+p'(t)[dx(t)]+[dp(t)]'[dx(t)] \\ \notag
=&\bigg\{\Big\{\dot{P}(t)x(t)+P(t)\Big[Ax(t)+\bar{A}\mathbb{E}[x(t)]+Bu(t-h)+\bar{B}\mathbb{E}[u(t-h)]\Big] \\ \notag
&+\dot{\bar{P}}(t)\mathbb{E}[x(t)] +\bar{P}(t)\Big[(A+\bar{A})\mathbb{E}[x(t)]+(B+\bar{B})\mathbb{E}[u(t-h)]\Big]+\Theta_1(t)\Big\}'x(t) \\ \notag
&+\Big[P(t)x(t)+\bar{P}(t)\mathbb{E}[x(t)]+\Theta(t)\Big]'\Big[Ax(t)+\bar{A}\mathbb{E}[x(t)]+Bu(t-h)+\bar{B}\mathbb{E}[u(t-h)]\Big] \\ \notag
&+\Big[Cx(t)+\bar{C}\mathbb{E}[x(t)]+Du(t-h)+\bar{D}\mathbb{E}[u(t-h)]\Big]' P(t) \\ \notag
&\times\Big[Cx(t)+\bar{C}\mathbb{E}[x(t)]+Du(t-h)+\bar{D}\mathbb{E}[u(t-h)]\Big]\bigg\}dt \\ \notag
&+\bigg\{\Big[Cx(t)+\bar{C}\mathbb{E}[x(t)]+Du(t-h)+\bar{D}\mathbb{E}[u(t-h)]\Big]'P(t)x(t)\\ \notag
&+p'(t)\Big[Cx(t)+\bar{C}\mathbb{E}[x(t)] +Du(t-h)+\bar{D}\mathbb{E}[u(t-h)]\Big]\bigg\}dw(t).
\end{align}
Taking expectation, and then taking integral from $h$ to $T$ on both sides of \eqref{26} and using Remark \ref{r1}, we have
\begin{align}\label{27}
&\mathbb{E}[p'(T)x(T)-p'(h)x(h)] \\ \notag
&=\mathbb{E}\int_h^T\bigg\{x'(t)\big[\dot{P}(t)+P(t)A+A'P(t)+C'P(t)C-\Pi(t,t+h)\big]x(t) \\ \notag
&+\mathbb{E}[x'(t)]\Big[\dot{\bar{P}}(t)+P(t)\bar{A} +\bar{A}'P(t)+\bar{P}(t)(A+\bar{A})+(A+\bar{A})'\bar{P}(t)+\bar{C}'P(t)C\\\notag
&+C'P(t)\bar{C}+\bar{C}'P(t)\bar{C}\Big]\mathbb{E}[x(t)] +2x'(t)\big[P(t)B+C'P(t)D\big]u(t-h) \\ \notag
&+2\mathbb{E}[x'(t)]\big[P(t)\bar{B}+\bar{P}(t)(B+\bar{B})+\bar{C}'P(t)D +C'P(t)\bar{D}+\bar{C}'P(t)\bar{D}\big]\mathbb{E}[u(t-h)]\\ \notag
&+x'(t)\bigg[\Pi(t,t)\hat{x}(t|t-h)-A'\Theta(t)\\ \notag
&-\int_t^{t+h}\Pi(t,\theta)d\theta\Big(\bar{A}\mathbb{E}[x(t)]+Bu(t-h)+\bar{B}\mathbb{E}[u(t-h)]\Big)\bigg] \\ \notag
&+\Big[Ax(t)+\bar{A}\mathbb{E}[x(t)]+Bu(t-h)+\bar{B}\mathbb{E}[u(t-h)]\Big]'\Theta(t) + u'(t-h)D'P(t)Du(t-h) \\ \notag
&+\mathbb{E}[u'(t-h)]\Big[\bar{D}'P(t)D+D'P(t)\bar{D}+\bar{D}'P(t)\bar{D}\Big]\mathbb{E}[u(t-h)]\bigg\}dt \\ \notag
&=\mathbb{E}\int_h^T\bigg\{-x'(t)Qx(t)-\mathbb{E}[x'(t)]\bar{Q}\mathbb{E}[x(t)]-u'(t-h)Ru(t-h) \\ \notag
&-\mathbb{E}[u'(t-h)]\bar{R}\mathbb{E}[u(t-h)] -\mathbb{E}[x'(t)]\Pi(t,t)Ex(t)+\mathbb{E}[x'(t)]\bar{\Pi}(t,t)\mathbb{E}[x(t)] \\ \notag
&+x'(t)\Pi(t,t)\hat{x}(t|t-h)+2x'(t)M'_1(t)u(t-h) \\ \notag
&+2\mathbb{E}[x'(t)]\big[M_2(t)-M_1(t)\big]'\mathbb{E}[u(t-h)]+u'(t-h)\Upsilon_1(t)u(t-h) \\ \notag
&+\mathbb{E}[u'(t-h)]\big[\Upsilon_2(t)-\Upsilon_1(t)\big]\mathbb{E}[u(t-h)]\bigg\}dt.
\end{align}
Since $\Upsilon_1(t)$ and $\Upsilon_2(t)>0$, we have
\begin{align}\label{28}
J_T=&\mathbb{E}\bigg\{\int_0^hx'(t)Qx(t)+\mathbb{E}[x'(t)]\bar{Q}\mathbb{E}[x(t)]dt+x'(h)p(h) \\ \notag
&+\int_h^T\bigg[u(t-h)-\mathbb{E}[u(t-h)]+[\Upsilon_1(t)]^{-1} M_1(t)\big[\hat{x}(t|t-h)-\mathbb{E}[x(t)]\big]\bigg]'\Upsilon_1(t) \\ \notag
&\times\bigg[u(t-h)-\mathbb{E}[u(t-h)]+[\Upsilon_1(t)]^{-1}M_1(t)\big[\hat{x}(t|t-h)-\mathbb{E}[x(t)]\big]\bigg] \\ \notag
&+\bigg[\mathbb{E}[u(t-h)] +[\Upsilon_2(t)]^{-1}M_2(t)\mathbb{E}[x(t)]\bigg]'\Upsilon_2(t) \\ \notag
&\times\bigg[\mathbb{E}[u(t-h)]+[\Upsilon_2(t)]^{-1}M_2(t)\mathbb{E}[x(t)]\bigg]dt\bigg\},
\end{align}
where $\Upsilon_1(t),~\Upsilon_2(t)$ and $M_1(t),~M_2(t)$ are given by \eqref{05}-(11) and \eqref{07}-(15).
From the above equation and the positive definiteness of $\Upsilon_1(t)$ and $\Upsilon_2(t)$, we can conclude that the optimal cost functional can be denoted as \eqref{091}, and the optimal controller satisfies
\begin{align}\label{29a}
&u(t-h)-\mathbb{E}[u(t-h)]+\big[\Upsilon_1(t)]^{-1}M_1(t)[\hat{x}(t|t-h)-\mathbb{E}[x(t)]\big]=0, \\ \label{29b}
&\mathbb{E}[u(t-h)]+[\Upsilon_2(t)]^{-1}M_2(t)\mathbb{E}[x(t)]=0.
\end{align}
The optimal controller \eqref{08} can be easily derived from \eqref{29a}-\eqref{29b}. This means that the sufficiency proof of the theorem is completed.

``Necessary":
We will prove that $\Upsilon_1(t)$ and $\Upsilon_2(t)>0, t\in[h,T]$ if the Problem \ref{p01} has a unique solution under Assumption \ref{a01}.
In fact, when Problem \ref{p01} is solvable, for purpose of deriving analytical form of $u(t)$ from the equilibrium equation \eqref{03}, we should attempt to find the relationship among $p(t)$, $q(t)$ and $x(t)$ in the first place. In fact, combining the earlier work of \cite{ZX} and \cite{Yong}, it is assume that
\begin{align}\label{11}
p(t)=P(t)x(t)+\bar{P}(t)\mathbb{E}[x(t)]-\int_0^h\mathbf{\Pi}(t,t+\theta)\hat{x}(t|t+\theta-h)d\theta,
\end{align}
where
\begin{align}\label{33}
-\dot{P}(t)=&Q+P(t)A+A'P(t)+C'P(t)C-\mathbf{\Pi}(t,t+h), \\
-\dot{\bar{P}}(t)=&\bar{Q}+P(t)\bar{A}+\bar{A}'P(t)+\bar{P}(t)(A+\bar{A})+(A+\bar{A})'\bar{P}(t)+C'P(t)\bar{C}\\\notag&+\bar{C}'P(t)C+\bar{C}'P(t)\bar{C}-\bar{A}'\int_0^h\mathbf{\Pi}(t,t+\theta)d\theta-\int_0^h\mathbf{\Pi}(t,t+\theta)d\theta\bar{A}\\\notag
&+\mathbf{\Pi}(t,t)-\mathbf{\bar{\Pi}}(t,t), \end{align}
and $\mathbf{\Pi}(t,t+h)$ and $\mathbf{\bar{\Pi}}(t,t)$ are given as below
\begin{align}\label{34}
&\mathbf{\Pi}(t,s)=e^{A'(s-t)}\Pi(s,s)e^{A(s-t)}, \quad s\in[t,t+h], \\
&\mathbf{\Pi}(t,t)=M'_1(t)\Upsilon_1(t)^{\dag}M_1(t), \\
&\mathbf{\bar{\Pi}}(t,t)=M'_2(t)\Upsilon_2(t)^{\dag}M_2(t),
\end{align}
with the terminal value $P(T)$ and $\bar{P}(T)$ as those in \eqref{02}.

Firstly, we show that $\Upsilon_1(t)$ and $\Upsilon_2(t)\geq 0$. Similar to the methods in \cite{QZ}. Suppose that the positive semidefiniteness does not hold, then there exist $t_1,~t_2\in[h,T]$ such that $\lambda_1(t_1)<0$, or $\lambda_2(t_2)<0$, where $\lambda_i(t)$ is a arbitrary eigenvalue of $\Upsilon_i(t)$, and the corresponding eigenvector is $v_i(t)~(i=1,2)$.

We choose $u_1(t-h)=L(t)\big[\hat{x}(t|t-h)-\mathbb{E}[x(t)]\big]+(L(t)+\bar{L}(t))\mathbb{E}[x(t)]$, where
\begin{align}\label{35}
&L(t)=\left\{\begin{array}{lll}\mathbf{K}(t), &\textrm{$\lambda_1(t) = 0,$}\\
\frac{\delta I_l^1}{|\lambda_1(t)|^{\frac{1}{2}}}v_1(t)+\mathbf{K}(t), &\textrm{$\lambda_1(t)\neq 0,$}\end{array}\right. \\
\\ \notag
&L(t)+\bar{L}(t)=\left\{\begin{array}{lll}\mathbf{K}(t)+\mathbf{\bar{K}}(t), &\textrm{$\lambda_2(t)=0,$}\\
\frac{\delta I_l^2}{|\lambda_2(t)|^{\frac{1}{2}}}v_2(t)+\mathbf{K}(t)+\mathbf{\bar{K}}(t), &\textrm{$\lambda_2(t)\neq 0,$}\end{array}\right.\\\notag
\end{align}
where $I_l^i$ is the indicator function of set $\{t\in[0,T]|\lambda_i(t)<-\frac{1}{l}\},l=1,2,\cdots, ~i=1,2$,  and $\textbf{K}(t),~\bar{\textbf{K}}(t)$ satisfy
\begin{align}\label{r001}
\textbf{K}(t)=&-\Upsilon_1(t)^{\dag}M_1(t),\\
\bar{\textbf{K}}(t)=&-\big[\Upsilon_2(t)^{\dag}M_2(t)-\Upsilon_1(t)^{\dag}M_1(t)\big].
\end{align}

Substituting controller $u_1(t)$ into \eqref{28}, we set the corresponding cost functional $J_T^1$ as:
\begin{align}\label{36}
J_T^1=&\mathbb{E}\bigg\{\int_0^hx'(t)Qx(t)+\mathbb{E}[x'(t)]\bar{Q}\mathbb{E}[x(t)]dt+x'(h)p(h) \\ \notag
&-\delta^2\int_{I_l^1}\Big(\hat{x}(t|t-h)-\mathbb{E}[x(t)]\Big)'\Big(\hat{x}(t|t-h)-\mathbb{E}[x(t)]\Big)dt \\ \notag
&-\delta^2\int_{I_l^2}\mathbb{E}[x'(t)]\mathbb{E}[x(t)]dt\bigg\}.
\end{align}
Assume that $\hat{x}(t|t-h)-\mathbb{E}[x(t)]\neq 0$ in the set $\{t\in[0,T]|\lambda_1(t)<-\frac{1}{l}\}$, $\mathbb{E}[x(t)]\neq 0$ in the set $\{t\in[0,T]|\lambda_2(t)<-\frac{1}{l}\}$.
Using the assumption that $\lambda_1(t_1)<0$, or $\lambda_2(t_2)<0$ and the continuity of the eigenvalues, we have $\mathcal{M}[I_l^1]+\mathcal{M}[I_l^2]>0$, where $\mathcal{M}$ is the Lebesgue measure. Therefore, $\int_{I_l^1}[\hat{x}(t|t-h)-\mathbb{E}[x(t)]]'[\hat{x}(t|t-h)-E[x(t)]]dt+\int_{I_l^2}\mathbb{E}[x'(t)]\mathbb{E}[x(t)]dt>0$. If we let $\delta\rightarrow +\infty$, we have $J_T^1\rightarrow -\infty$, which is in contradiction with $J_T^1\geq 0$.
Therefore, we have $\Upsilon_1(t)$ and $\Upsilon_2(t)\geq 0$.

Next, we show that $\Upsilon_1(t)$ and $\Upsilon_2(t)>0$. Since the condition does not hold, it follows from $\Upsilon_1(t)$, $\Upsilon_2(t)\geq 0$ and \eqref{28} that we choose
\begin{align*}
 u^{(2)}(t)=\mathbf{K}(t)x(t)+\mathbf{\bar{K}}(t)Ex(t),~~~u^{(3)}(t)=\mathbf{K}(t)x(t)+\mathbf{\bar{K}}(t)Ex(t)+\mathbf{\bar{L}},
\end{align*}
where $\bar{\mathbf{L}}=\big[I-\Upsilon_2(t)^\dag\Upsilon_2(t)\big]\bar{z}$, and $\bar{z}\neq 0$, thus $ u^{(2)}(t)\neq  u^{(3)}(t)$. Substituting $u^{(2)}(t)$ and $u^{(3)}(t)$ into \eqref{28} respectively, the corresponding cost functionals satisfy
\begin{align}\label{77}
J_T^2=J_T^3=\mathbb{E}\bigg\{\int_0^hx'(t)Qx(t)+\mathbb{E}[x'(t)]\bar{Q}\mathbb{E}[x(t)]dt+x'(h)p(h)\bigg\}=J_T^*.
\end{align}
Since Problem \ref{p01} is uniquely solvable, we obtain $u^{(2)}(t)=u^{(3)}(t)$, which is a contradiction to $u^{(2)}(t) \neq u^{(3)}(t)$. Therefore, we can conclude that $\Upsilon_1(t)$ and $\Upsilon_2(t)>0$.
\qed

\subsection{Appendix D: Proof of Theorem \ref{t03}}
\textit{Proof}.
The outline of the proof can be summarized as follows:
\begin{itemize}
  \item $P_T(t)$ is monotonically increasing with respect to $T$.
  \item $P_T(t)$ is uniformly bounded.
  \item $\bar{P}_T(t)$ and $\displaystyle\int_0^h\mathbf{\Pi}_T(t,t+\theta)d\theta$ are convergent with respect to $T$.
  \item $P_T(t)$ and $\bar{P}_T(t)$ satisfy
$\lim\limits_{T\rightarrow \infty}\dot{P}_T(t)=0 \mbox{ and }
\lim\limits_{T\rightarrow \infty}\dot{\bar{P}}_T(t)=0.$
  \item $\mathbf{\Pi}_T(t,t+\theta),\theta\in[0,h]$ and $\bar{\mathbf{\Pi}}_T(t,t)$ are convergent with respect to $T$.
  \item $P-\int_0^h\mathbf{\Pi}(\theta)d\theta > 0$ and $P+\bar{P}-\displaystyle\int_0^h\mathbf{\Pi}(\theta)d\theta > 0$.
\end{itemize}

Firstly, we show that $P_T(t)$ is monotonically increasing with respect to $T$. For any initial state $x_0,~\mathbb{E}[x_0]=0$, and $u(\tau)=0,~\tau\in[-h,0)$, Denote the corresponding optimal cost functional as $J_T^0$. Taking the integral from $0$ to $T$ on both side of \eqref{27}, with $P(T)$ and $\bar{P}(T)=0$, we have
\begin{align}\label{43}
J_T^0= ~ &\mathbb{E}\bigg\{[P_T(0)x(0)+\bar{P}_T(0)\mathbb{E}[x(0)]-\int_0^h\mathbf{\Pi}_T(0,\theta)\hat{x}(0|\theta-h)d\theta]'x(0) \\ \notag
&+\int_0^h\Big[\hat{x}(t|t-h)-\mathbb{E}[x(t)]\Big]'M'_1(t)\Upsilon_1(t)^{\dag}M_1(t)\Big[\hat{x}(t|t-h)-\mathbb{E}[x(t)]\Big]dt \\ \notag
& +\int_0^h\mathbb{E}[x'(t)]M'_2(t)\Upsilon_2(t)^{\dag}M_2(t)\mathbb{E}[x(t)]dt\bigg\}.
\end{align}

For $t\in[0,h]$, we have $d\mathbb{E}[x(t)]=(A+\bar{A})\mathbb{E}[x(t)]dt$. Therefore, $\mathbb{E}[x(t)]=e^{(A+\bar{A})t}\mathbb{E}[x(0)]=0$. Using \eqref{17}, we get $\hat{x}(t|s)=\hat{x}(t|0)=e^{At}x_0$. Hence,
\begin{align}\label{44}
J_T^0=&\mathbb{E}\bigg\{x'_0\bigg[P_T(0)x_0-\int_0^h\mathbf{\Pi}_T(0,\theta)d\theta\bigg]x_0+\int_0^hx'_0e^{A't}M'_1(t)\Upsilon_1(t)^{\dag}M_1(t)e^{At}x_0dt\bigg\} \\ \notag
=&\mathbb{E}[x'_0P_T(0)x_0]\geq 0.
\end{align}
In fact, for any $t<T_1<T_2$, we have
\begin{align}\label{45}
\mathbb{E}[x'_0P_{T_1}(0)x_0]=J_{T_1}^0\leq J_{T_2}^0=\mathbb{E}[x'_0P_{T_2}(0)x_0].
\end{align}
Since $x_0$ is arbitrary, $0 \leq P_{T_1}(0)\leq P_{T_2}(0)$. Using time-invariance of $P_T(t)$, we obtain that for any $t<T_1<T_2$,
$P_{T_1}(t)=P_{T_1-t}(0)\leq P_{T_2-t}(0)= P_{T_2}(t)$, i.e., $P_T(t)$ monotonically increases with respect to $T$ and $P_T(t)\geq 0$.

Next, we show that $P_T(t)$ is uniformly bounded. For this purpose, we rewrite system \eqref{01} as the following system:
\begin{align}\label{47}
d\mathbb{X}(t)=[\mathbb{A}\mathbb{X}(t)+\mathbb{B}\mathbb{U}(t-h)]dt+[\mathbb{C}\mathbb{X}(t)+\mathbb{D}\mathbb{U}(t-h)]dw(t),
\end{align}
where 
$$\begin{array}{lll}
\mathbb{X}(t)=(x(t),\mathbb{E}[x(t)])', & \mathbb{U}(t)=(u(t),\mathbb{E}[u(t)])',
& \mathbb{A}=\left[
  \begin{array}{ccc}
    A&\bar{A}\\
    0&A+\bar{A}\\
  \end{array}
\right], \\ [5mm]
\mathbb{B}=\left[
  \begin{array}{ccc}
    B&0\\
    0&B+\bar{B}\\
  \end{array}
\right], & 
\mathbb{C}=\left[
  \begin{array}{ccc}
    C&C+\bar{C}\\
    0&0\\
  \end{array}
\right], &
\mathbb{D}=\left[
  \begin{array}{ccc}
    D&D+\bar{D}\\
    0&0\\
  \end{array}
\right].
\end{array}$$
The mean-square stabilizability of system \eqref{01} with controller $u(t-h)=K\hat{x}(t|t-h)+\bar{K}\mathbb{E}[x(t)]$ equals to the stabilization of system \eqref{47} with controller
$$\mathbb{U}(t-h)=\left(
  \begin{array}{ccc}
    K&\bar{K} \\
    0&K+\bar{K} \\
    \end{array}
\right)
\left(
  \begin{array}{ccc}
    \hat{x}(t|t-h) \\
    \mathbb{E}[x(t)] \\
    \end{array}
\right)=\mathbb{K}\left(
  \begin{array}{ccc}
    \hat{x}(t|t-h) \\
    \mathbb{E}[x(t)] \\
    \end{array}
\right).$$
 Similar to the method for the stochastic delay systems in \cite{ZX}, an equivalent abstract model for system \eqref{47} without delay can be represented as as follows. Let $f(t)=(\mathbb{X}(t),\mathbb{U}(t,s))$, where $\mathbb{U}(t,s)=\mathbb{U}(t+s),~s\in[-h,0]$. For any initial state $x_0,~\mathbb{E}[x_0]=0$, and $u(\tau)=0,~\tau\in[-h,0)$, system \eqref{47} can be changed into the following abstract model:
\begin{align}\label{46}
\left\{\begin{array}{lll}df(t)=[\mathcal{A}f(t)+\mathcal{B}\mathbb{U}(t)]dt+[\mathcal{C}f(t)+\mathcal{D}\mathbb{U}(t)]dw(t), \\
f(0)=(\mathbb{X}(0),0),
\end{array}\right.
\end{align}
where $\mathcal{A}$ is the infinitesimal generator of a strongly continuous semigroup over $\mathcal{M}_2$, defined by:
\begin{align}\label{48}
\mathcal{A}\eta=\mathcal{A}(\eta_0,\eta_1)=\bigg(\mathbb{A}\eta_0+\mathbb{B}\eta_1(-h),\frac{d}{ds}\eta_1\bigg),
\end{align}
where the domain $\mathcal{D}(\mathcal{A})=\{(\eta_0,\eta_1)\in \mathcal{M}_2:\eta_1(\theta)\in W^{1,2}(-h,0;\mathbb{R}^m),\eta_1(0)=0\}$.
$u(t)$ is a control input taking values in a real separable Hilbert space $\mathcal{U}$, and the remaining operators are defined as:
\begin{align}\label{49}
&\mathcal{B}\mathbb{U}(t)=(0,\delta_0(s)\mathbb{U}(t)),\\
&\mathcal{C}f(t)=(\mathbb{C}\mathbb{X}(t)+\mathbb{D}\mathbb{U}(t-h),0),\\
&\mathcal{D}\mathbb{U}(t)=(0,0).
\end{align}

Given the stabilizability of system \eqref{47}, the infinite-dimensional abstract system \eqref{46} is stabilizable, i.e., there exists a stabilizing controller $\mathbb{U}(t)=\mathcal{K}f(t)$ such that $\mathbb{E}\int_0^\infty\langle f(t),f(t)\rangle dt<\infty$. Therefore, there exists $c_0$ such that
\begin{align}\label{50}
\int_0^\infty \mathbb{E}[x'(t)]\mathbb{E}[x(t)]dt   \leq
\int_0^\infty \mathbb{E}[x'(t)x(t)]dt  \leq
c_0 x'_0x_0, \\ \label{50a}
\int_h^\infty \mathbb{E}[\hat{x}'(t|t-h)\hat{x}(t|t-h)]dt   \leq
\int_0^\infty \mathbb{E}[x'(t)x(t)]dt  \leq
c_0 x'_0x_0.
\end{align}
Combining the above two inequalities with \eqref{37}, we have
\begin{align}\label{51}
J_T^0 &\leq J=\mathbb{E}\bigg\{\int_0^\infty x'(t)Qx(t)+\mathbb{E}[x'(t)]\bar{Q}\mathbb{E}[x(t)]+u'(t)Ru(t)+\mathbb{E}[u'(t)]\bar{R}\mathbb{E}[u(t)]dt\bigg\} \\ \notag
&=\mathbb{E}\bigg\{\int_0^\infty x'(t)Qx(t)+\mathbb{E}[x'(t)]\bar{Q}\mathbb{E}[x(t)]dt \\ \notag
& \quad +\int_h^\infty \big[\hat{x}(t|t-h)-\mathbb{E}[x(t)]\big]'K'RK\big[\hat{x}(t|t-h)-\mathbb{E}[x(t)]\big]   \\ \notag
&\quad +\mathbb{E}[x'(t)](K+\bar{K})'(R+\bar{R})(K+\bar{K})\mathbb{E}[x(t)]dt\bigg\} \\ \notag
&\leq \mathbb{E}\bigg\{\int_0^\infty x'(t)Qx(t)+\mathbb{E}[x'(t)]\bar{Q}\mathbb{E}[x(t)]
+\hat{x'}(t|t-h)K'RK \hat{x}(t|t-h) \\ \notag
& \quad +\mathbb{E}[x'(t)](K+\bar{K})'(R+\bar{R})(K+\bar{K})\mathbb{E}[x(t)]dt\bigg\} \\ \notag
&\leq \lambda_{max}\big\{Q+\bar{Q}+K'RK+(K+\bar{K})'(R+\bar{R})(K+\bar{K})\big\}c_0 x'_0x_0.
\end{align}
Therefore, with \eqref{44} and the arbitrariness of $x_0$, the convergence of $P_T(t)$ is proved, i.e. there exists $P\geq0$ such that $\displaystyle\lim_{T\rightarrow \infty}P_T(t)=P$.

Furthermore, we show the convergence of $\bar{P}_T(t)$ and $\int_0^h\mathbf{\Pi}_T(t,t+\theta)d\theta$.
In fact, let us consider system
\eqref{01} starting at $t=h$ with an arbitrary initial value $x(h)$, and the corresponding cost functional is given by
\begin{align}\label{29}
J_T(h)=&\mathbb{E}\bigg\{\int_h^T x'(t)Qx(t)+\mathbb{E}[x'(t)]\bar{Q}\mathbb{E}[x(t)]+ u'(t-h)Ru(t-h)\\\notag
&+\mathbb{E}[u'(t-h)]\bar{R}\mathbb{E}[u(t-h)]dt +x'(T)P(T)x(T)+\mathbb{E}[x'(T)]\bar{P}(T)\mathbb{E}[x(T)]\bigg\}.
\end{align}
Similarly, the optimal value of cost functional \eqref{29} is as follows
\begin{align}\label{30}
J^*_T(h)=&\mathbb{E}\bigg\{x'(h)\bigg[P_T(h)x(h)+\bar{P}_T(h)\mathbb{E}[x(h)]-\int_0^h\Pi_T(h,h+\theta)\hat{x}(h|\theta)d\theta\bigg]\bigg\} \\ \notag
=&\mathbb{E}\bigg\{x'(h)\bigg[P_T(h)-\int_0^h\Pi_T(h,h+\theta)d\theta\bigg]x(h)+\mathbb{E}[x'(h)]\bar{P}_T(h)\mathbb{E}[x(h)]\bigg\}.
\end{align}
If we choose an arbitrary $x(h)$, $\mathbb{E}[x(h)]=x(h)$, we have
\begin{align}\label{31}
\mathbb{E}\bigg\{x'(h)\bigg[P_T(h)+\bar{P}_T(h)-\int_0^h\Pi_T(h,h+\theta)d\theta\bigg]x(h)\bigg\}\geq0.
\end{align}
In addition, if we choose an arbitrary $x(h)$, $\mathbb{E}[x(h)]=0$, we have
\begin{align}\label{32}
\mathbb{E}\bigg\{x'(h)\bigg[P_T(h)-\int_0^h\Pi_T(h,h+\theta)d\theta\bigg]x(h)\bigg\}\geq0.
\end{align}
Therefore, given that $x(h)$ is arbitrary,  we get $P_T(h)+\bar{P}_T(h)-\int_0^h\Pi_T(h,h+\theta)d\theta\geq0$ and $P_T(h)-\int_0^h\Pi_T(h,h+\theta)d\theta\geq0$. Similar to the discussions above, we know that $P_T(h)+\bar{P}_T(h)-\int_0^h\Pi_T(h,h+\theta)d\theta$ and $P_T(h)-\int_0^h\Pi_T(h,h+\theta)d\theta$ increase with respect to $T$. By the same method as the second step, \eqref{50}-\eqref{50a} hold under the initial values $\mathbb{E}[x(h)]=x(h)$ and $\mathbb{E}[x(h)]=0$, respectively. This implies the boundness. Therefore, the convergence of $\bar{P}_T(t)$ and $\int_0^h\mathbf{\Pi}_T(t,t+\theta)d\theta$ can be proved similarly.

Next, we prove that $P_T(t)$ and $\bar{P}_T(t)$ satisfy
\begin{equation}\label{52}
\lim_{T\rightarrow \infty}\dot{P}_T(t)=0 \mbox{ and }
\lim_{T\rightarrow \infty}\dot{\bar{P}}_T(t)=0.
\end{equation}
In fact, similar to the discussions in \cite{ZX}, it is easy to know that $\dot{P}_T(t)$ and $\dot{\bar{P}}_T(t)$ are uniformly continuous. Applying Barbalat's Lemma (Lemma 8.2 in \cite{Khalil}), we get the desired result \eqref{52}.

Then, we prove that $\mathbf{\Pi}_T(t,t+\theta),\theta\in[0,h]$ and $\bar{\mathbf{\Pi}}_T(t,t)$ are convergent with respect to $T$. In fact, taking $T \rightarrow \infty$ on both side of \eqref{33}, together with \eqref{52}, we know that $\displaystyle\lim_{T\rightarrow\infty}\mathbf{\Pi}_T(t,t+h)$ exists. We denote the limitation as $\mathbf{\Pi}(h)$.

Since
\begin{align}\label{53}
\mathbf{\Pi}_T(t,t+\theta)&=e^{A'\theta}\mathbf{\Pi}_T(t+\theta,t+\theta)e^{A\theta}\\\notag
&=e^{A'\theta}\mathbf{\Pi}_{T+h-\theta}(t+h,t+h)e^{A\theta}\\\notag
&=e^{-A'(h-\theta)}\mathbf{\Pi}_{T+h-\theta}(t,t+h)e^{-A(h-\theta)},
\end{align}
we have $\displaystyle\lim_{T\rightarrow\infty}\mathbf{\Pi}_T(t,t+\theta)=\mathbf{\Pi}(\theta)$.
Obviously, $\displaystyle\lim_{T\rightarrow\infty}\bar{\mathbf{\Pi}}_T(t,t)=\bar{\mathbf{\Pi}}(0)$ holds.

Furthermore, we will show that $P-\int_0^h\mathbf{\Pi}(\theta)d\theta > 0$ and $P+\bar{P}-\int_0^h\mathbf{\Pi}(\theta)d\theta > 0$.
From \eqref{31}-\eqref{32} and the arbitrariness of $x(h)$, it is obvious that they are positive semidefinite. Thus, if $P-\int_0^h\mathbf{\Pi}(\theta)d\theta$ is not positive definite, there must exist $x_0\neq 0$ and $\mathbb{E}[x_0]=0$ such that
\begin{align}\label{54}
\mathbb{E}\bigg(x'_0\bigg[P-\int_0^h\mathbf{\Pi}(\theta)\hat{x}(t|t+\theta-h)d\theta \bigg)'x_0\bigg]=0.
\end{align}
Using It\^{o}'s formula  to $\Big[Px(t)+\bar{P}\mathbb{E}[x(t)]-\int_0^h \mathbf{\Pi}(\theta)\hat{x}(t|t+\theta-h)d\theta\Big]'x(t)$ and then taking expectation, we have
\begin{align}\label{55}
&\mathbb{E}\bigg\{d \bigg[Px(t)+\bar{P}\mathbb{E}[x(t)]-\int_0^h \mathbf{\Pi}(\theta)\hat{x}(t|t+\theta-h)d\theta\bigg]'x(t) \bigg\} \\ \notag
=&\mathbb{E}\big\{-x'(t)Qx(t)-\mathbb{E}[x'(t)]\bar{Q}\mathbb{E}[x(t)]-u'(t-h)Ru(t-h) \\ \notag
&-\mathbb{E}[u'(t-h)]\bar{R}\mathbb{E}[u(t-h]-\mathbb{E}[x'(t)]\Pi(0)\mathbb{E}[x(t)]+\mathbb{E}[x'(t)]\bar{\Pi}(0)\mathbb{E}[x(t)] \\ \notag
&+x'(t)\Pi(0)\hat{x}(t|t-h)+2x'(t)M'_1u(t-h)+2\mathbb{E}[x'(t)]\big[M_2-M_1\big]'\mathbb{E}[u(t-h] \\ \notag
&+u'(t-h)\Upsilon_1u(t-h)+\mathbb{E}[u'(t-h)]\big[\Upsilon_2-\Upsilon_1\big]\mathbb{E}[u(t-h)]\big\},
\end{align}
where $M_1,~M_2,~\Upsilon_1,~\Upsilon_2$ are presented by \eqref{40}-\eqref{40c}. Taking integration on both sides of the above equation form $0$ to $T$,
\begin{align}\label{56}
&\mathbb{E}\bigg\{\int_0^T x'(t)Qx(t)+\mathbb{E}[x'(t)]\bar{Q}\mathbb{E}[x(t)]dt+\int_h^T u'(t-h)Ru(t-h)\\\notag
&+\mathbb{E}[u'(t-h)]\bar{R}\mathbb{E}[u(t-h)]dt \bigg\} \\ \notag
=&-\mathbb{E}\bigg\{\bigg[Px(T)+\bar{P}\mathbb{E}[x(T)]-\int_0^h \mathbf{\Pi}(\theta)\hat{x}(T|T+\theta-h)d\theta\bigg]'x(T)\bigg\} \\ \notag
&+\mathbb{E}\bigg\{\bigg[Px(0)+\bar{P}\mathbb{E}[x(0)]-\int_0^h \mathbf{\Pi}(\theta)x(0)d\theta\bigg]'x(0)\bigg\} \\ \notag
&-\mathbb{E}\bigg\{\int_0^h u'(t-h)Ru(t-h)+\mathbb{E}[u'(t-h)]\bar{R}\mathbb{E}[u(t-h)]dt \bigg\} \\ \notag
&+\mathbb{E}\bigg\{\int_0^T \bigg[u(t-h)-\mathbb{E}[u(t-h)]+[\Upsilon_1]^{-1}M_1\big[\hat{x}(t|t-h)-\mathbb{E}[x(t)]\big]\bigg]' \\ \notag
&\times\Upsilon_1\bigg[u(t-h)-\mathbb{E}[u(t-h)] +[\Upsilon_1]^{-1}M_1\big[\hat{x}(t|t-h)-\mathbb{E}[x(t)]\big]\bigg] \\ \notag
&+\bigg[\mathbb{E}[u(t-h)]+[\Upsilon_2]^{-1}M_2\mathbb{E}[x(t)]\bigg]'\Upsilon_2\bigg[\mathbb{E}[u(t-h)]+[\Upsilon_2]^{-1}M_2\mathbb{E}[x(t)]\bigg]dt\bigg\}.
\end{align}
We choose the initial condition as $x(0)=x_0,~u(\tau-h)=K\hat{x}(\tau|\tau-h)+\bar{K}\mathbb{E}[x(\tau)],\tau\in[0,h)$, and the stabilizing controller as $u(t-h)=K\hat{x}(t|t-h)+\bar{K}\mathbb{E}[x(t)],t\in[h,T]$, together with \eqref{54}, the above equation can be rewritten as
\begin{align}\label{57}
&\mathbb{E}\bigg\{\int_0^T x'(t)Qx(t)+\mathbb{E}[x'(t)]\bar{Q}\mathbb{E}[x(t)]dt+\int_h^T u'(t-h)Ru(t-h) \\ \notag
&+\mathbb{E}[u'(t-h)]\bar{R}\mathbb{E}[u(t-h)]dt\bigg\}\\\notag
=&-\mathbb{E}\bigg\{\bigg[Px(T)+\bar{P}\mathbb{E}[x(T)]-\int_0^h \mathbf{\Pi}(\theta)\hat{x}(T|T+\theta-h)d\theta\bigg]'x(T) \bigg\} \\ \notag
&-\mathbb{E}\bigg\{\int_0^h u'(t-h)Ru(t-h) +\mathbb{E}[u'(t-h)]\bar{R}\mathbb{E}[u(t-h)]dt\bigg\} \\ \notag
=&-\mathbb{E}\bigg\{\big[x(T)-\mathbb{E}[x(T)]\big]'\bigg[P-\int_0^h\mathbf{\Pi}(\theta)d\theta\bigg]\big[x(T)-\mathbb{E}[x(T)]\big] \\ \notag
&+\mathbb{E}[x'(T)]\bigg[P+\bar{P}-\int_0^h\mathbf{\Pi}(\theta)d\theta\bigg]\mathbb{E}[x(T)]\bigg\} \\ \notag
&-\mathbb{E}\bigg\{\int_0^h\big[u(t-h)-\mathbb{E}[u(t-h)]\big]'R\big[u(t-h)-\mathbb{E}[u(t-h)]\big] \\ \notag
&+\mathbb{E}[u'(t-h)][R+\bar{R}]\mathbb{E}[u(t-h)]dt\bigg\} \leq 0,
\end{align}
where
\begin{align*}
\mathbb{E}\bigg\{x'(T)\int_0^h \mathbf{\Pi}(\theta)\mathbb{E}[x(T)|\mathcal{F}_{T+\theta-h}]d\theta\bigg\}
=&\mathbb{E}\bigg\{\mathbb{E}\bigg[x'(T)\int_0^h \mathbf{\Pi}(\theta)d\theta x(T)|\mathcal{F}_{T+\theta-h}\bigg]\bigg\} \\ \notag
=&\mathbb{E}\bigg\{x'(T)\int_0^h \mathbf{\Pi}(\theta)d\theta x(T)\bigg\}
\end{align*}
Therefore, it follows from \eqref{57} that we have
\begin{align}\label{r002}
0 \leq &\mathbb{E}\bigg\{\int_0^T x'(t)Qx(t)+\mathbb{E}[x'(t)]\bar{Q}\mathbb{E}[x(t)]dt+\int_h^T u'(t-h)Ru(t-h)\\ \notag
&+\mathbb{E}[u'(t-h)]\bar{R}\mathbb{E}[u(t-h)]dt \bigg\} \\ \notag
=&\mathbb{E}\bigg\{\int_0^T \big[x(t)-\mathbb{E}[x(t)]\big]'Q\big[x(t)-\mathbb{E}[x(t)]\big]+\mathbb{E}[x'(t)]\big(Q+\bar{Q}\big)\mathbb{E}[x(t)]dt\\\notag
&+\int_h^T \big[u(t-h)-\mathbb{E}[u(t-h)]\big]'R\big[u(t-h)-\mathbb{E}[u(t-h)]\big] \\ \notag
&+\mathbb{E}[u'(t-h)]\big(R+\bar{R}\big)\mathbb{E}[u(t-h)]dt\bigg\} \leq 0,
\end{align}
which indicates $\mathcal{Q}^{\frac{1}{2}}X(t)=0,a.s.,\forall t\in[0,T]$. From the exact observability of system $(A,~\bar{A},~C,~\bar{C},~\mathcal{Q}^{\frac{1}{2}})$, we have
$X(0)=\left(
  \begin{array}{ccc}
    x_0-\mathbb{E}[x_0] \\
    \mathbb{E}[x_0] \\
  \end{array}
\right)=0$,
which is a contradiction to $x_0\neq 0$. Thus, $P-\int_0^h\mathbf{\Pi}(\theta)d\theta > 0$. Similarly, the positive definiteness of $P+\bar{P}-\int_0^h\mathbf{\Pi}(\theta)d\theta$ can be achieved. Due to space limitation, the proof is omitted here. Thereby, the proof of Theorem \ref{t03} is completed.
\qed

\subsection{Appendix E: Proof of Theorem \ref{t04}}
\textit{Proof}.
"Necessity":
Suppose that system \eqref{01} is stabilizable in the mean-quare sense, according to Theorem \ref{t03}, the existence of a positive definite solution to the coupled algebraic Riccati equations \eqref{38}-\eqref{38b} is established. Then we only need to prove the uniqueness.
In fact, assume that there are two different solutions $(P_1, \bar{P}_1)$ and $(P_2, \bar{P}_2)$, where $P_1>0,~P_1+\bar{P}_1>0$ and $P_2>0,~P_2+\bar{P}_2>0$, and the associated feedback gains are $K_1,~\bar{K}_1$ and $K_2,~\bar{K}_2$. In particular, if we choose the initial values as $x(0)=x_0, ~ \mathbb{E}[x_0]=0$, $u(t-h)=0,~t\in[0,h)$, in this case, the optimal cost functional is
\begin{align}\label{62}
J^*=&\mathbb{E}\bigg\{x'(0)\bigg[P_1x(0)+\bar{P}_1\mathbb{E}[x(0)]-\int_0^h \mathbf{\Pi}_1(\theta)d\theta x(0)\bigg] \\ \notag
&+\int_0^h \hat{x}'(t|t-h)\mathbf{\Pi}_1(0)\times\hat{x}(t|t-h)dt\bigg\}\\\notag
=&\mathbb{E}[x'_0P_1x_0].
\end{align}
Similarly, we have $J^*=\mathbb{E}[x'_0P_2x_0]$. By the arbitrariness of $x_0$, we can conclude that $P_1=P_2$. $P$ is uniquely solvable from \eqref{38}, which implies that $\mathbf{\Pi}(\theta)$ can be uniquely obtained by \eqref{39}-\eqref{39b}.
When system \eqref{01} starts from $t=h$, the cost functional can be presented by
\begin{align}\label{63}
J(h)=&\mathbb{E}\bigg\{\int_h^\infty x'(t)Qx(t)+\mathbb{E}[x'(t)]\bar{Q}\mathbb{E}[x(t)]+ u'(t-h)Ru(t-h)\\\notag
&+\mathbb{E}[u'(t-h)]\bar{R}\mathbb{E}[u(t-h)]dt\bigg\}.
\end{align}
Following the method from \eqref{29} to \eqref{31} for this case, the optimal cost functional reduces to
\begin{align}\label{64}
J^*(h)=\mathbb{E}\bigg\{x'(h)\bigg[P_1+\bar{P}_1-\int_0^h\mathbf{\Pi}_1(\theta)d\theta\bigg]x(h)\bigg\}.
\end{align}
Similarly, we conclude that $P+\bar{P}-\int_0^h\mathbf{\Pi}(\theta)d\theta$ is unique due to the arbitrariness of $x(h)$. Thus $\bar{P}$ is unique.

"Sufficiency":
In this part, we show that the stochastic mean-field system with delay is stabilizable in the mean-square sense with the specific controller \eqref{59}. It is equivalent to show that the following system
\begin{align}\label{66}
dx(t)=&[Ax(t)+BK\hat{x}(t|t-h)+(\bar{A}+B\bar{K}+\bar{B}K+\bar{B}\bar{K})\mathbb{E}[x(t)]]dt\\\notag
&+[Cx(t)+DK\hat{x}(t|t-h)+(\bar{C}+D\bar{K}+\bar{D}K+\bar{D}\bar{K})\mathbb{E}[x(t)]]dw(t)
\end{align}
is mean-square stable.

Firstly, the Lyapunov function candidate is defined as follows:
\begin{align}\label{65}
V(t,x(t))=\mathbb{E}\bigg[x'(t)Px(t)+\mathbb{E}[x'(t)]\bar{P}\mathbb{E}[x(t)]-x'(t)\int_0^h\mathbf{\Pi}(\theta)\hat{x}(t|t+\theta-h)d\theta\bigg],~~t\geq h.
\end{align}
If the coupled algebraic Riccati equations have a positive definite solution $(P, \bar{P})$, we get
\begin{align}\label{67}
V(t,x(t)) &\geq \mathbb{E}\bigg\{x'(t)Px(t)+\mathbb{E}[x'(t)]\bar{P}\mathbb{E}[x(t)]-x'(t)\int_0^h\mathbf{\Pi}(\theta)\hat{x}(t|t+\theta-h)d\theta \\ \notag
& \quad -\int_0^h\tilde{x}'(t|t+\theta-h)\mathbf{\Pi}(\theta)\tilde{x}(t|t+\theta-h)d\theta\bigg\} \\ \notag
&=\mathbb{E}\bigg\{x'(t)Px(t)+\mathbb{E}[x'(t)]\bar{P}\mathbb{E}[x(t)]-x'(t)\int_0^h\mathbf{\Pi}(\theta)d\theta x(t)\bigg\} \\ \notag
&=\mathbb{E}\bigg\{\big[x(t)-\mathbb{E}[x(t)]\big]'\bigg[P-\int_0^h\mathbf{\Pi}(\theta)d\theta\bigg]\big[x(t)-\mathbb{E}[x(t)]\big]+\mathbb{E}[x'(t)]\\\notag
&~~~\times\bigg[P+\bar{P}-\int_0^h\mathbf{\Pi}(\theta)d\theta\bigg]\mathbb{E}[x(t)]\bigg\} \geq 0,
\end{align}
where $x(t)=\hat{x}(t|t+\theta-h)+\tilde{x}(t|t+\theta-h)$, i.e.,
$\tilde{x}(t|t+\theta-h)$ denotes the error between the state $x(t)$ and the estimation $\hat{x}(t|t+\theta-h)$.
 Taking derivation to $V(t,x(t))$, we have
\begin{align}\label{68}
&dV(t,x(t)) \\ \notag
= & \mathbb{E}\Big\{-x'(t)Qx(t)-\mathbb{E}[x'(t)]\bar{Q}\mathbb{E}[x(t)] \\ \notag
&-u'(t-h)Ru(t-h) -\mathbb{E}[u'(t-h)]\bar{R}\mathbb{E}[u(t-h)] \\ \notag
&+\big[u(t-h)-K\hat{x}(t|t-h)-\bar{K}\mathbb{E}[x(t)]\big]'\Upsilon_1\big[u(t-h)-K\hat{x}(t|t-h)-\bar{K}\mathbb{E}[x(t)]\big] \\ \notag
&+\big[\mathbb{E}[u(t-h)]-(K+\bar{K})\mathbb{E}[x(t)]\big]'\Upsilon_2\big[\mathbb{E}[u(t-h)]-(K+\bar{K})\mathbb{E}[x(t)]\big]\Big\}dt \\ \notag
=&\mathbb{E}\Big\{-x'(t)Qx(t)-\mathbb{E}[x'(t)]\bar{Q}\mathbb{E}[x(t)]-u'(t-h)Ru(t-h) \\ \notag
&-\mathbb{E}[u'(t-h)]\bar{R}\mathbb{E}[u(t-h)]\Big\}dt \leq 0.
\end{align}
Combining \eqref{67}-\eqref{68}, we know that $V(t,x(t))$ decreases and is bounded. Hence, $V(t,x(t))$ is convergent when $t\rightarrow\infty$.

Taking integration on both sides of \eqref{68}, we have
\begin{align}\label{69}
&V(h,x)-V(T,x) \\ \notag
=&\mathbb{E}\bigg\{\int_h^T -x'(t)Qx(t)-\mathbb{E}[x'(t)]\bar{Q}\mathbb{E}[x(t)]-u'(t-h)Ru(t-h) \\ \notag
&-\mathbb{E}[u'(t-h)]\bar{R}\mathbb{E}[u(t-h)]dt\bigg\} \\ \notag
\geq & \mathbb{E}\bigg\{\int_h^T -{x^*}'(t)Qx^*(t)-\mathbb{E}[{x^*}'(t)]\bar{Q}\mathbb{E}[x^*(t)]-{u^*}'(t-h)Ru^*(t-h)\\ \notag&-\mathbb{E}[{u^*}'(t-h)]\bar{R}\mathbb{E}[u^*(t-h)]dt\bigg\},
\end{align}
where $u^*(t-h)=K(t)\hat{x}^*(t|t-h)+\bar{K}(t)\mathbb{E}[x^*(t)]$ is the optimal controller for the finite horizon problem (from Theorem \ref{t02}), and $x^*(t)$ is the corresponding state trajectory. From \eqref{29}-\eqref{30}, \eqref{69} can be rewritten as
\begin{align}\label{70}
V(h,x)-V(T,x)\geq & \mathbb{E}\bigg\{x'(h)P_T(h)x(h)+\mathbb{E}[x'(h)]\bar{P}_T(h)\mathbb{E}[x(h)] \\ \notag
&-x'(h)\int_0^h\Pi_T(h,h+\theta)\hat{x}(h|\theta)d\theta\bigg\}\geq 0.
\end{align}
Similarly, we have
\begin{align}\label{71}
&V(t+h,x)-V(t+T,x)\\\notag
\geq & \mathbb{E}\bigg\{x'(t+h)P_{T}(h)x(t+h)+\mathbb{E}[x'(t+h)]\bar{P}_T(h)\mathbb{E}[x(t+h)] \\ \notag
&-x'(t+h)\int_0^h\Pi_T(h,h+\theta)\hat{x}(t+h|t+\theta)d\theta\bigg\}\\\notag
\geq & \mathbb{E}\bigg\{\big[x(t+h)-\mathbb{E}[x(t+h)]\big]'\Big[P_T(h)-\int_0^h\Pi_T(h,h+\theta)d\theta\Big]\big[x(t+h)-\mathbb{E}[x(t+h)]\big] \\ \notag
&+\mathbb{E}\big[x'(t+h)\big]\Big[P_T(h)+\bar{P}_T(h)-\int_0^h\Pi_T(h,h+\theta)d\theta\Big]\mathbb{E}\big[x(t+h)\big]\bigg\}\geq 0.
\end{align}
Next, we show that $P_T(h)-\int_0^h\Pi_T(h,h+\theta)d\theta>0$. If it is not this case, then there must exist $x_h,~x_h\neq 0$, $\mathbb{E}[x_h]=0$ such that $\mathbb{E}\bigg\{x'_h\bigg[P_T(h)-\int_0^h\Pi_T(h,h+\theta)d\theta\bigg]x_h\bigg\}=0$. If we
choose $h$ be the initial time and $x_h$ be the initial state, we have
\begin{align}\label{72}
&\mathbb{E}\bigg\{\int_h^T {x^*}'(t)Qx^*(t)+\mathbb{E}[{x^*}'(t)]\bar{Q}\mathbb{E}[x^*(t)]+{u^*}'(t-h)Ru^*(t-h)\\\notag
&+\mathbb{E}[{u^*}'(t-h)]\bar{R}\mathbb{E}[u^*(t-h)]dt\bigg\} \\ \notag
=&\mathbb{E}\bigg\{x'_h\bigg[P_T(h)-\int_0^h\Pi_T(h,h+\theta)d\theta\bigg]x_h\bigg\}=0.
\end{align}
From the above equation, we can see that $\mathcal{Q}^{\frac{1}{2}}X(t)=0, a.s., \forall t\in[h,T]$. From Lemma \ref{l2} and the exact observability of $(A,~\bar{A},~C,~\bar{C},~\mathcal{Q}^{\frac{1}{2}})$, we can conclude that $x_h=0$, which is a contradiction. Therefore, $P_T(h)-\int_0^h\Pi_T(h,h+\theta)d\theta>0$. Similarly, we can prove that $P_T(h)+\bar{P}_T(h)-\int_0^h\Pi_T(h,h+\theta)d\theta>0$.

Thus, \eqref{71} can be further formulated as
\begin{align}\label{73}
&V(t+h,x)-V(t+T,x)\\\notag
\geq &\lambda_0 \mathbb{E}\{[x(t+h)-\mathbb{E}[x(t+h]]'[x(t+h)-\mathbb{E}[x(t+h)]]+\mathbb{E}[x'(t+h)x(t+h)]\} \\ \notag
=& \lambda_0 \mathbb{E}[x'(t+h)x(t+h)]\geq 0,
\end{align}
where
$\lambda_0  =\min\big\{\lambda(P_T(h)-\int_0^h\Pi_T(h,h+\theta)d\theta>0),\lambda(P_T(h)+\bar{P}_T(h)-\int_0^h\Pi_T(h,h+\theta)d\theta)\big\}>0$. Taking $t\rightarrow \infty$ on both sides of the above equation leads to
$\displaystyle\lim_{t\rightarrow \infty}\mathbb{E}[x'(t)x(t)]=0$. This implies the mean-square stability of system \eqref{01}.

In the following, we present the optimal controller for infinite horizon problem and the corresponding optimal cost functional. Similar to \eqref{56}, applying It\^{o}'s formula to $\bigg[Px(t)+\bar{P}\mathbb{E}[x(t)]-\int_0^h \mathbf{\Pi}(\theta)\hat{x}(t|t+\theta-h)d\theta\bigg]'x(t)$, we get
\begin{align}\label{74}
&\mathbb{E}\bigg\{\int_0^T x'(t)Qx(t)+\mathbb{E}[x'(t)]\bar{Q}\mathbb{E}[x(t)]dt \\ \notag
&+\int_h^T u'(t-h)Ru(t-h)+\mathbb{E}[u'(t-h)]\bar{R}\mathbb{E}[u(t-h)]dt\bigg\} \\ \notag
=&\mathbb{E}\bigg\{\bigg[Px(0)+\bar{P}\mathbb{E}[x(0)]- \int_0^h\mathbf{\Pi}(\theta)x(0)d\theta\bigg]'x(0) \\ \notag
&-\bigg[Px(T)+\bar{P}\mathbb{E}[x(T)]-\int_0^h\mathbf{\Pi}(\theta)\hat{x}(T|T+\theta-h)d\theta\bigg]'x(T) \\\notag
&+\int_0^h\Big\{-\mathbb{E}[x'(t)]\Pi(0)\mathbb{E}[x(t)]+\mathbb{E}[x'(t)]\bar{\Pi}(0)\mathbb{E}[x(t)]
+\hat{x}'(t|t-h)\Pi(0)\hat{x}(t|t-h) \\ \notag
&+2x'(t)M'_1u(t-h)+2\mathbb{E}[x'(t)][M_2-M_1]'\mathbb{E}[u(t-h)]+u'(t-h)\Upsilon_1u(t-h) \\ \notag
&+\mathbb{E}[u'(t-h)][\Upsilon_2-\Upsilon_1]\mathbb{E}[u(t-h)]\Big\}dt \\ \notag
&+\int_h^T \bigg[u(t-h)-\mathbb{E}[u(t-h)]+[\Upsilon_1]^{-1}M_1\big[\hat{x}(t|t-h)-\mathbb{E}[x(t)]\big]\bigg]' \Upsilon_1 \\ \notag
&\times\bigg[u(t-h)-\mathbb{E}[u(t-h)] +[\Upsilon_1]^{-1}M_1\big[\hat{x}(t|t-h)-\mathbb{E}[x(t)]\big]\bigg] \\ \notag
& +\bigg[\mathbb{E}[u(t-h)]+[\Upsilon_2]^{-1}M_2\mathbb{E}[x(t)]\bigg]'\Upsilon_2\bigg[\mathbb{E}[u(t-h)]+[\Upsilon_2]^{-1}M_2\mathbb{E}[x(t)]\bigg]dt\bigg\}.
\end{align}
Taking $T\rightarrow\infty$ in the above equation, together with $\displaystyle\lim_{t\rightarrow\infty}\mathbb{E}[x'(t)x(t)]=0$, we get
\begin{align}\label{75}
J=&\mathbb{E}\bigg\{\bigg[Px(0)+\bar{P}\mathbb{E}[x(0)]-\int_0^h \mathbf{\Pi}(\theta)x(0)d\theta\bigg]'x(0) \\ \notag
&+\int_0^h\Big\{-\mathbb{E}[x'(t)]\Pi(0)\mathbb{E}[x(t)]+\mathbb{E}[x'(t)]\bar{\Pi}(0)\mathbb{E}[x(t)] \\ \notag
&+\hat{x}'(t|t-h)\Pi(0)\hat{x}(t|t-h)+2x'(t)M'_1u(t-h) \\ \notag
&+2\mathbb{E}[x'(t)][M_2-M_1]'\mathbb{E}[u(t-h)] +u'(t-h)\Upsilon_1u(t-h) \\ \notag
&+\mathbb{E}[u'(t-h)]]\Upsilon_2-\Upsilon_1]\mathbb{E}[u(t-h)]\Big\}dt \\ \notag
&+\int_h^\infty \bigg[u(t-h)-\mathbb{E}[u(t-h)]+[\Upsilon_1]^{-1}M_1\big[\hat{x}(t|t-h)-\mathbb{E}[x(t)]\big]\bigg]'\Upsilon_1 \\ \notag
&\times\bigg[u(t-h)-\mathbb{E}[u(t-h)] +[\Upsilon_1]^{-1}M_1\big[\hat{x}(t|t-h)-\mathbb{E}[x(t)]\big]\bigg] \\ \notag
&+\bigg[\mathbb{E}[u(t-h)]+[\Upsilon_2]^{-1}M_2\mathbb{E}[x(t)]\bigg]'\Upsilon_2\bigg[\mathbb{E}[u(t-h)]+[\Upsilon_2]^{-1}M_2\mathbb{E}[x(t)]\bigg]dt\bigg\}.
\end{align}
Hence, it follows from \eqref{75} that we have the desired optimal controller \eqref{59} and its corresponding optimal cost functional \eqref{61}.
\qed

\end{document}